\magnification=\magstep1
\input amstex

\documentstyle{amsppt}
\pageheight{8.5truein}
\NoBlackBoxes
\leftheadtext{Monotonic Local Decay Estimates}
\rightheadtext{Soffer}

\topmatter
\title  Monotonic Local Decay Estimates
\endtitle
\author  Avy Soffer\endauthor
\address Mathematics Department,
 Rutgers University, 110 Frelinghuysen road, Piscataway , NJ 08854
\endaddress
\email soffer at math.rutgers.edu
\endemail
\abstract  For the Hamiltonian operator $H=-\Delta + V(x)$ of the Schr\"odinger Equation with a
repulsive potential, the problem of local decay is considered. It is analyzed by a direct method, based
on a new, $L^2$ bounded, propagation observable. The resulting decay estimate,
 is in certain cases {\bf monotonic} in time, with no ``Quantum
Corrections''.  This method is then applied to some examples in one
and higher dimensions.  In particular the case of the Wave Equation on
a Schwarzschild manifold is redone: Local decay, stronger than the
known ones are proved (minimal loss of angular derivatives and lower order of
radial derivatives of initial data).  The method developed here can be an alternative
in some cases to the Morawetz type estimates, with
$L^2$-multipliers replacing the first order operators.  It provides an
alternative to Mourre's method, by including thresholds and high energies.
\endabstract
\endtopmatter

\document

\head Section 1\endhead
\subhead
1. Introduction\endsubhead
\medskip

The starting point to a-priori estimates for dispersive equations is
finding an operator which generates a monotonic function relative to
the flow; the prime examples are  Morawetz identity, the Dilation
identity and the pseudo conformal identity.  The Morawetz identity
applies in three or more dimensions.

The above identities are generated by differential operators $M,$ and we
have

$$
\frac{d}{dt}\langle \psi(t), M\psi(t)\rangle \geq 0\tag1.1
$$
where $\psi(t)$ is the solution of Schr\"odinger Equation at time $t$, the $ \langle,\rangle$ stands for the
usual $L^2$-scalar product.

To derive the Morawetz estimate we choose ($n$-dimension)
$$
M =- i \frac{\partial}{\partial r} - \frac{n-2}{r}, n \geq 3, \qquad r\equiv |x|.
\tag 1.2
$$
The Dilation identity: $M= \frac{1}{2} (x\cdot p+p\cdot x)$
$$
p=-i\nabla_x\tag 1.3
$$
The Conformal identity
$$
-M = |\frac{x}{t} - p|^2 t^\alpha +t^{2-\alpha}V(x) \quad,  0\leq \alpha \leq 2. \tag 1.4
$$
The aim of this note is to construct {\bf monotonic} observables $M$ which
are microlocal or phase-space operators.

The implications of such a construction include new local
decay estimates, in particular, in one dimension, and new propagation
estimates; it opens the way to new  classes of  a-priori
estimates, including local decay at thresholds.

The operators $M$ which I refer to, for obvious reason, as {\bf propagation observables} (PROB)
, are also known as {\bf multipliers}.

\head Section 2 \endhead
\subhead
2a. Some notation and preliminaries  \endsubhead
\medskip

We consider the Schr\"odinger flow on $L^2(\Bbb R^n)$ generated
by a self-adjoint operator $H$:
$$
i\frac{\partial \psi}{\partial t} = H \psi ; \qquad \psi (t=0) \in L^2.\tag 2.1
$$

We will focus on the case where
$$
\aligned
H = & -\Delta + V (x)\\
    & -\Delta  \qquad \text{ is the Laplacian on } \Bbb R^n.
\endaligned\tag 2.2
$$
We assume from now on that $V(x)$ is a real valued, uniformly
bounded $C^1$ function of $x\in \Bbb R^n$, so that $H$ is self-adjoint on the domain
$D(-\Delta ) = H^2(\Bbb R^n),$ the Sobolev space.

In $  L^2(\Bbb R^n)$, we define the momentum operator $p,$

$$
p=-i\partial_x \text{ and } r=|x|.
$$
Then,
$$
-\Delta = p\cdot p \equiv p^2.\tag 2.4
$$
We let
$$
A= \frac{1}{2}(x\cdot p +p \cdot x) = x\cdot p-n i /2 = p \cdot x +
n i /2 = \left(-ir \frac{\partial}{\partial r} - i
\frac{\partial}{\partial r} r \right) / 2 \tag 2.5
$$
and we have that
$$
\imath[p_i, x_j] = \delta_{ij}\, ; \qquad i,j=1,...,n ,\tag 2.6
$$
where $\delta_{ij}$ stands for the Kronecker delta function.
Therefore:
$$
\aligned
&  \imath[-\Delta, x_j]= 2p_j\,  ;\quad  \imath[A, x_j] = x_j \\
&  \imath [A, p_j]= - p_j\, ; \quad
\imath [-\Delta, A] = 2
p^2.\endaligned
\tag 2.7
$$

We denote $\langle x\rangle^2 = 1+|x|^2$ and by $F(B\in I)$ the
smoothed projection of the self-adjoint operator $B$ in the interval
$I$.  E.g.,
$$
F(|x|\leq 1) \text{ stands for the multiplier}
$$
by the smoothed characteristic function of $I\equiv \{|\lambda | \leq 1\}$ in
$L^2(\Bbb R^n_x)$.

From equation (2.7) we derive ,
$$
e^{\imath sA}p_je^{-\imath sA}=e^{-s}p_j \, ;\quad  e^{\imath sA}x_je^{-\imath sA}=e^sx_j.
$$
\subhead
2b- Monotonic propagation Estimates \endsubhead
 It is generally known, from the works of Enss and Mourre  that scattering states propagate into becoming "outgoing".
 So, in particular, one can prove, using the Mourre estimate that
 $$
 \|P^-(A)e^{-iHt}g(H)P_c(H)\psi\| \leq o(t),
 $$
 as $t$ approaches $+\infty.$
 Here, $P^-(A)$ is the projection on the negative spectral part of the Dilation generator, $A.$
 When $g(H)\equiv 1,$ we get decay with essentially no rate.
 When $g(H)$ is supported away from zero and infinity, one can prove fast decay in time, for localized initial data in space, as well as minimal and maximal velocity bounds [HSS and cited ref].
 It is much more difficult to get estimates when the cutoff $g$ is not present, and no localization of the initial data is assumed.
 In this case the methods of Mourre and [HSS] do not apply, in general.
 Some generalizations were obtained in [Ger, MRT and cited ref.], see also [Rod-T], replacing the Mourre estimate with a weak version of it.
 Here, I will develop a new way of getting decay estimates, for certain classes of hamiltonians, without localizing $H$ or $\psi.$
 Furthermore, I will show that the propagation from the region of incoming waves into outgoing waves, and similar propagation estimates, is monotonic in time, for the free flow, and for the free flow perturbed by a class of repulsive potentials. These are two typical results:
  I show that it is possible to modify, by exponentially small corrections at infinity, the projection $P^{-}(A)$ so that, the solution decays monotonically on its range, for the repulsive potentials:
 $$
\langle \psi(t), F^-_M(A) \psi (t) \rangle \downarrow 0,\text{ as } t\to
+ \infty
$$
and
$$
\int^T_0 \|\langle x\rangle^{-1}F(A\leq -M) \psi (t) \|^2 dt \leq \langle \psi(0), 2 F^-_M(A)
\psi(0)\rangle,
$$
 see Proposition (6.2).
 The first part shows that the flow from incoming waves to outgoing is monotonic, with no restriction on the initial data!
 The second estimate shows, that at least locally in space, the incoming part is controlled, integrably in time, by the size of the incoming waves part of the initial data. So, in particular, no incoming wave can reappear locally,
 including zero energy and high energy contributions.
 The above estimates hold in any dimension, including one dimension, for one hump potentials.
 This has immediate applications to the case of scattering of the wave equation on Black-hole metrics:
 \medskip
 \proclaim{Theorem }
For the Hamiltonian with the Schwarzschild potential $\ell^2 V(x),$ with $V$ analytic repulsive, we have the following estimate:
$$
\int_0^T \|F(|x|\leq r_0)\ell u\|^2 \leq c\ln \ell E(u).
$$
\endproclaim
See section 8.
Previously, a similar estimate was obtained in[B-Sof3,4], by complicated multi-step phase space propagation estimates.
The propagation estimates above extends to time dependent hamitonians, with small, sufficiently localized potential perturbations.
\head Section 3 \endhead
\subhead
3. The propagation observable  \endsubhead
\medskip

Since $A,$ the dilation generator defined in equation (2.5), is a self-adjoint operator, we can construct the operator
$F(A/R)$:
$$
F\left(\frac{A}{R}\right) \equiv \tanh \frac{A}{R}
\tag 3.1
$$
by the spectral theorem.

We show that $F(A/R)$ has a positive commutator with $H=-\Delta$, and
find lower bounds for it, if $R$ is sufficiently large.

Then, this is extended to $H=-\Delta + V$ for certain classes of
potentials $V$.

Note that the analysis works in any dimension, and we specify to one
dimension, which is the more difficult case.

To proceed, recall the commutator expansion Lemma [Sig-Sof1-2].

Let
$$
ad_A^n(B)\equiv [ad_A^{n-1}(B),A] ;\quad ad_A^1=[B,A].
$$

\proclaim{Lemma 3.1 } Commutator Expansion Lemma

$$
\aligned
& i[B, f(A)] = \int \hat f(\lambda)e^{i\lambda A}[e^{-i\lambda A}Be^{i\lambda A}-B] d\lambda\\
&= f'(A) i [B, A] + \frac{1}{2!} f'' (A) i [[B, A], A] +
\cdots R_n
\endaligned\tag 3.2
$$
$$
R_n = \frac{1}{n!}\int \hat f (\lambda) e^{i\lambda A}
\int^{\lambda}_0 e^{-isA} \int^s_0  e^{-i\mu A}\cdots\int^t_0 e^{-iu
A}(-\imath)^n ad^n_B (A) e^{+iu A} du \ldots d\lambda. \tag 3.3
$$
\endproclaim

In particular, we get:
 \proclaim{Corollary 3.2}

Let $A$ be the dilation generator as defined before, on $L^2(\Bbb R^n).$
$$
\text{For }\quad  R>2/ \pi
$$
$$
\tanh A/R : D(-\Delta) \to D(-\Delta).
$$
\endproclaim
\demo{Proof}
Commuting $\Delta $ through $e^{i\lambda A/R}$, we have:
$$
 e^{i\lambda A/R}[\Delta ,e^{-i\lambda A/R}]=e^{i\lambda A/R}\Delta e^{-i\lambda A/R} - \Delta = (e^{-2\lambda
/R} - 1 ) \Delta : D(\Delta) \to L^2.
$$
Therefore, using the Commutator Expansion Lemma
with $n=1,$ and the property (3.6) of the Fourier Transform of the tanh function, the result follows.
\enddemo

\proclaim{Theorem 3.3}  $i[-\Delta, \tanh (A/R)] = 2 pg^2 (A/R) p
\geq 0$, for $R>2/\pi$.
Here,
$$
g^2(A/R)=\frac{\sin (2/R)}{{\cosh \frac{2A}{R}+ 2\cosh \frac{2}{R}}}.
$$
\endproclaim

\demo{Proof}
$$
\text{In the sense of forms, on} \quad D(H) \times D(H):
$$
$$
\aligned
&i[p^2, \tanh (A/R)] = i p [p,\tanh (A/R)] + i [p, \tanh (A/R)] p\\
&= p \int \hat f (\lambda) e^{i\lambda A/R} (-i) \int^{\lambda}_0 e^{-i
s A/R}(p/R) e^{i s A/R} ds d\lambda + c.c.  \\
&= p\int \hat f (\lambda) e^{i\lambda A/R}(-i) \int^{\lambda}_0
e^{+s/R} (p/R) dsd \lambda + c.c.\\
&=-i p\int \hat f (\lambda) e^{i\lambda A/R} p \,  e^{+s/R}
\big|^{\lambda}_0 d\lambda + c.c.\\
&=-i p\int \hat f (\lambda) \left(e^{\lambda/R}-1\right) e^{i\lambda
A/R} p \, d \lambda  + c.c.\\
&= -ip\left[\tanh \left( \frac{A+i}{R}\right) - \tanh \left(\frac{A}{R}\right)\right] p + c.c.\\
&=p\frac{1}{i} \left[\tanh \frac{A+i}{R} - \tanh \frac{A-i}{R} \right] p
\endaligned \tag 3.4
$$
provided $|\hat f(\lambda) | \leq c e^{-k|\lambda|}$ with $ k >\frac
{1}{R}, |\lambda| > 1.$

We also note that
$$
\hat f (\lambda) \left(e^{\lambda/R} -1\right) \sim \frac{1}{\lambda}
\left(e^{\lambda/R} -1\right)\text{ near zero, }\tag 3.5
$$
which is bounded.

$$
\hat f(\lambda) =  \frac{\pi}{\sinh \pi \lambda} , \lambda > 0,\tag
3.6
$$
and similar formula for $ \lambda <0.$
$$
\aligned &\frac{1}{i} \left( \tanh \frac{A + i}{R} - \tanh
\frac{A-i}{R} \right) = \frac{1}{i} {\frac{\sinh(2i/R)}{\cosh
\frac{A+ i}{R} \cosh
\frac{A-i}{R}}}\\
&= \frac{\sin (2/R)}{{\cosh \frac{2A}{R}+ 2\cosh \frac{2}{R}}} > 0 \quad
\text{ for } R>2/\pi.\endaligned
\tag 3.7
$$
\enddemo

\qed

\proclaim{Corollary 3.4}{Propagation estimate}

 For $R>2/\pi, H =- \Delta
$
$$
\aligned
&\left\langle \psi (t) ,  \, \tanh \frac{A}{R} \psi (t) \right\rangle - \left\langle\psi
(0), \tanh \frac{A}{R} \psi (0)\right\rangle\\&
 = \int^t_0 d s \| g(A) p \psi (s)
   \|^2 \leq 2\|\psi\|_{L^2}^2.
\endaligned \tag 3.8
$$
with $g^2(A) \geq \frac{C}{R} \frac{1}{\cosh\frac{2A}{R}}.$
\endproclaim
\demo{proof}
$$
\text{For }\quad \psi \in D(H):
$$
$$
\aligned &\frac{d}{dt}\left\langle e^{-iHt}\psi,\tanh \frac{A}{R}
e^{-iHt}\psi \right\rangle \\
 &= \left\langle H\psi(t),\tanh\frac{A}{R}\psi(t)\right\rangle - \left\langle
\psi(t),\tanh\frac{A}{R}H\psi(t)\right\rangle \\
&=\left\langle \psi(t),\left(H\tanh
\frac{A}{R}-\tanh\frac{A}{R}H\right)\psi(t)\right\rangle.
\endaligned
$$
The first equality follows by Von Neumann's Theorem. The second
equality follows by Corollary 3.2 and Spectral Theorem. The
Corollary now follows from Theorem 3.3 and Fundamental Theorem of
Calculus.
\enddemo
Few remarks are in order now.
\medskip

\noindent{\bf Remark 1.}  The above estimate shows that in the region
$|A|\leq C$, the solution has an extra gain of {\bf one} derivative,
upon time averaging.  One expects, more generally, that away from the
propagation set, in the phase-space, that the gain in derivatives
should be high.

Another important conclusion is the monotonicity of the flow in the
phase space.
\medskip

\noindent{\bf Remark 2.}  The corollary implies that the left hand
side is monotonically increasing in time, in fact, with non vanishing
derivative.

This means that the flow from the region $A\leq 0$ to the region
$A\geq 0$ is strictly monotonic.  This has important applications:

Define
$$
F^+_M (A/R) \equiv \left( F\left(\frac{A-M}{R}\right) + 1\right)
/2.\tag 3.9
$$
Then, the function $F^+_M(A/R)$ is exponentially close to the
projection operator $P^+ (A\geq M)$.
$
\left(\text{ for $|A|$ large enough depending on $R$}\right),
$
the projection on outgoing waves.

We can then immediately conclude that outgoing part of the solution
is strictly monotonic increasing up to exponentially small
correction of order $e^{-M}$.

Moreover, since the solution decays in time in the complement region,
we see that the correction is $o (t) e^{-M}$.

This property will remain true under decaying potential perturbations,
in some sense, since for large $M$, the potential term is
$O(M^{-\sigma})$ is this region.

\head Section 4 \endhead
\subhead
4. Adding Potentials \endsubhead
\medskip

The main interest in this note will be the case of ``one hump''
potentials in one dimension.  These are repulsive potentials V, such
that
$$
i[V, A] = - x\cdot \nabla V \geq 0.\tag 4.1
$$
We begin with the simple model
$$
V_0(x) = \frac{c_0}{b^2 + x^2}  ,\qquad   c_0 > 0.\tag 4.2
$$

Then, we have that Monotonic propagation estimates hold for $H_0 = -
\Delta + V_0 (x)$:

\proclaim{Proposition 4.1}
For $H_0=-\Delta +V_0(x),$ as above,
$$
i[H_0, F(A/R)] = 2p g^2 (A/R)] p + c_0 \frac{2}{b^2 + x^2} x \, g^2(A/R)
x \frac{1}{b^2 + x^2} \tag 4.3
$$
$$
\aligned
& 2\int^t_0 \|g(A/R) p\psi (s) \|^2 ds + c_0 \int^t_0 \|g(A/R)
\frac{x}{b^2 + x^2} \psi (s) \|^2 ds\\
& =\langle \psi (t) , \tanh (A/R) \psi (t) \rangle - \langle \psi (0),
\tanh (A/R) \psi (0)\rangle
\endaligned
\tag 4.4
$$
and $g^2(A/R)\gtrsim \frac{C}{R} \cosh^{-1} (2 A/R)$ , as before.
\endproclaim

\demo{Proof}
The proof follows from Theorem 3.3 and its application with $x$
replacing $p$:
$$
 ( \text{ Thm } 3.3 ) : \qquad  i[-\Delta, F(A/R)] = 2 p g^2 (A/R) p
$$
$$
\aligned
&
i[V_0, F(A/R)] = + c_0 \frac{2}{b^2 + x^2} i[F(A),  x^2]
\frac{1}{b^2 + x^2} \\
&= c_0 \frac{2}{b^2 + x^2} x g^2(A/R)x \frac{1}{b^2
+ x^2}
\endaligned
$$
where we use that $i[F(A), x^2] = 2x g^2(A) x$, the sign reversed
when $x\leftrightarrow p$. Equation (4.4) follows upon integrating over time
the Heisenberg identity for the Schr\"odinger equation.   \qed
\enddemo

The above theorem, and its proof, extends in a variety of situations:

\proclaim{Corollary 4.2}
Let
$$ H =- \Delta +V(x),
$$
and suppose that $V(x)$ admits a representation of the form:
$$
V(x) = \int^\infty_0 \frac{\rho(\alpha) d\alpha}{\alpha + x^2} ,\qquad \rho
(\alpha) \geq 0
$$
$\rho(\alpha)$ a positive measure, $|\rho(\alpha)|\leq c|\alpha|,
|\alpha|\leq 1$. We assume, moreover, that
$$
\int_0^{\infty} \frac{\rho(t)}{1+t} dt < \infty.
$$

  Then, the estimates of Theorem 4.1 hold for $H$,
with a different weight function in $x$:
$$
\frac{x}{b^2 + x^2} \to W_\rho(x)
$$
so that
$$
2\int^t_0 \|g(A/R) p\psi (s)\|^2 ds + c \int^t_0 \|g(A/R) W_\rho(x)
\psi (s) \|^2 ds
$$
$$
\leq \vert \langle \psi (t) ,  \, \tanh (A/R) \psi (t)\rangle \vert +\vert \langle\psi
(0), \tanh (A/R) \psi (0)\rangle \vert \tag 4.5
$$
\endproclaim

\noindent{\bf Remark } The class of potentials $V(x)$ above are Stieltjes functions.

\demo{Proof}
  The contribution from the potential term $V$ to the commutator is
  computed as before, to be
$$
\int^t_0 ds \int^{\infty}_0 \rho(\alpha)\| g(A/R)
  \frac{x}{\alpha + x^2} \psi(s)\|^2 d \alpha \tag 4.6
$$
$$
g(A/R) \frac{x}{\alpha + x^2} \psi = g \left(1-\frac{\alpha + x^2}{\alpha_0
+ x^2}\right) g^{-1} g \frac{x}{\alpha + x^2} \psi + g
\frac{x}{\alpha_0 + x^2} \psi \tag 4.7
$$
Now, if we integrate over $|\alpha-\alpha_0| < \delta
|\alpha_0|, \delta << 1$, we have that
$$
\aligned
& (\alpha_0 >0), \int_{|\alpha-\alpha_0|< \delta \alpha_0}
\|g(A/R)\frac{x}{\alpha + x^2} \psi \|^2 \rho(\alpha) d\alpha \geq \\
& \int_{|\alpha-\alpha_0|< \delta \alpha_0}\|g(A/R)\frac{x}{\alpha +
x^2} \psi \|^2 \rho(\alpha) d\alpha\\
&-c\int_{|\alpha-\alpha_0|< \delta \alpha_0}
\|g\frac{\alpha-\alpha_0}{\alpha_0 + x^2} g^{-1}\|^2
  \,  \|g\frac{x}{x^2 + \alpha^2}\psi\|^2 d\alpha.
\endaligned
\tag 4.8
$$
\enddemo
So, we only need to get smallness of
$$
\sup_{|\alpha-\alpha_0|<\delta \alpha_0} \|g
\frac{\alpha-\alpha_0}{\alpha + x^2} g^{-1} \| \leq
\sup_{|\alpha-\alpha_0|<\delta \alpha_0} 2 \delta
\|\frac{\alpha_0}{\alpha_0 + e^{2i\beta} x^2}\|_{L^\infty_x}
$$
 since in  our case $g^{-1} \sim  \frac{e^{-A/R} +  e^{A/R}}{2}$ and so $\beta \sim 1/R$.

So, for $R>1$, the result follows. Summing over the intervals
$\alpha$ around $\alpha_0= \frac{k}{N}$, for some large $N$, $0<k$
integer, we get a lower bound on the expression (4.6) of the form

$$
\aligned
\int^t_0 ds   &C \sum_{k} \int_{|\alpha_k - \frac{k}N|< \delta
k/N} \rho(\alpha)  \|g (A/R) \frac{x}{\alpha_k + x^2}
\psi (s) \|^2 d \alpha \\
=&C\int^t_0\sum_k \|g(A/R) \frac{\rho_k x}{\alpha_k + x^2} \psi (s)
\|^2  ds \\
\geq &C\int^t_0  ds \| g(A/R) \sum_k \frac{\rho_k / \langle k
\rangle^{1/2 + \varepsilon}x}{\alpha_k + x^2} \psi \|^2 \equiv
C\int^t_0 \|g(A/R) W_{\rho}(x) \psi (s) \|^2 ds.
\endaligned
\tag 4.9
$$
$$
\rho_k=\int_{|\alpha_k - \frac{k}N|< \delta
k/N} \rho(\alpha) d\alpha
$$

Next we need a {\bf microlocal uncertainty principle} inequality:

\proclaim{Lemma 4.3 }
For all $R$ large enough, $g$ a bounded $C^{\infty}$ function, $g(A/R)> 0,$
with,
$$
\sum_{i=1}^N \vert g^{(i)}\vert\leq c\vert g \vert,
$$
for sufficiently large $N=N(\sigma)>2,$ we have:

\medskip

\noindent
$$
(a) \qquad  (1+\varepsilon)  p g^2(A/R) p \geq g p^2 g \tag 4.10
$$

$
(b) \qquad   p^2 + \langle x \rangle^{-\sigma}_b x^2 \langle x \rangle^{-\sigma}_b\geq
$
$$
\frac{1}{2} \langle x \rangle^{-\sigma}_b (p^2 + x^2)\langle x
\rangle^{-\sigma}_b \geq \frac{1}{4} \langle x\rangle^{-2 \sigma}_b
\tag 4.11
$$
$ (c) \qquad pg^2(A/R) p + \langle x\rangle^{-\sigma}_b x g^2(A/R)
x\langle x\rangle^{-\sigma}_b \geq
\frac{1-\varepsilon}{4}g(A/R)\langle x\rangle ^{-2\sigma}_b g(A/R)$
 \newline for all $b$ large enough.
\endproclaim

\demo{Proof}
Part c) is proved using parts  a) b).  Assuming part b), we prove a)
and c):

$$
\aligned
& g = g(A/R); \qquad p = - i\nabla_x\\
& pg^2 p = pgp + [p,g] gp = gp^2 g + gp [g,p] + [p, g] gp\\
& =gp^2 g+ [gp,[g,p]] = gp^2 g +g [p, p\tilde g] + [g, p\tilde g^*] p\\
& = gp^2 g + gp\tilde g^{(2)*} p - p\tilde g \tilde g^* p.
\endaligned\tag 4.12
$$
$$
\tilde g\equiv [g,p], \qquad \tilde g^{(2)}\equiv [p,\tilde g].
$$

So, since by construction $\tilde g = O (\frac{1}{R}), \quad  \tilde g^{(2)} = O
\left(\frac{1}{R^2}\right)$, we get
$$
\aligned
p(g^2 + \tilde g \tilde g^*) p &= gp^2 g + gp \tilde g^{(2)*}p\\
&= gp^2 g + pg\tilde g^{(2) *} p + p\tilde g^ *\tilde g^{(2)*} p
\endaligned
$$
so,
$$
p(g^2 + \tilde g \tilde g^* - 2 Re g\tilde g^{(2)*}) p = g p^2 g.\tag
4.13
$$
Finally, for $R$ large,
$$
p(g^2 + \tilde g \tilde g^* - 2  Re  g\tilde g^{(2) *}) p \leq
(1+\varepsilon_R)p g^2 p\tag 4.14
$$
since $g>0$, vanishing only at infinity, and since $\tilde g, \tilde
g^{(2)}$ decay faster at $\infty$, and are of order $\frac{1}{R}$ and
$\frac{1}{R^2}$ respectively.  We therefore conclude that part a)
follows:

$$
(1+ \varepsilon_R) pg^2 p \geq gp^2 g. \tag 4.15
$$
\enddemo
Next, we prove part c):
\demo{Proof of c}
It follows from (4.15) that,
$$
\aligned
& pg^2(A/R) p + \langle x \rangle^{-\sigma}_b x g^2(A/R) x \langle
x\rangle^{-\sigma}_b \\
& \geq (1-\varepsilon) g(A/R) p^2 g(A/R) + (1 - \varepsilon)\langle
x\rangle^{-\sigma}_b g(A/R) x^2 g(A/R) \langle x\rangle^{-\sigma}_b.
\endaligned\tag 4.16
$$

We now need to commute $\langle x\rangle^{-\sigma}_b$ through
$g(A/R)$.  Commuting powers of $\langle x\rangle^{-1}_b$ through, the error
commutators terms are of the form

$\left(\psi, \langle x
\rangle^{-\sigma}_b\left(\frac{x}{\langle x \rangle_b}\right)^j\tilde g_1
a (x) x^2 \tilde g_2\langle x \rangle^{-\sigma}_b\left(
\frac{x}{\langle x\rangle_b}\right)^{j'}\psi \right), |a(x)|\leq c.$

Any such term is therefore bounded by
$$
\aligned
& c\left\{\| x \tilde g_1\left(\frac{x}{\langle x\rangle_b}\right)^j\langle
x \rangle^{-\sigma}_b \psi \|^2 + \| x \tilde g_2\left(
\frac{x}{\langle x\rangle}_b\right)^{j'} \langle
x \rangle^{-\sigma}_b \psi \|^2\right\}.
\endaligned
$$

For all $f(x)$, we have:
$$
\aligned
&\|x\tilde g_1f(x)\psi\| \leq \|{\tilde g}_1 xf (x) \psi \| + \| \tilde{\tilde
g}_1 x f(x)\psi \|\\
&\leq O\left(\frac{1}{R}\right) \|g x f(x) \psi \|,
\endaligned
$$
since $  \tilde g_1 = O\left(\frac{1}{R}\right),\quad \tilde{\tilde g}_1=O\left(\frac{1}{R^2}\right) ,$
and $ \vert g'\vert +\vert g''\vert \leq c\vert g\vert.$

Applying this last inequality with $f(x) = \left(\frac{x}{\langle x
\rangle_b} \right)^j \langle x \rangle^{-\sigma}_b$ we have that
$$
\|x\tilde g_1\left(\frac{x}{\langle x \rangle_b} \right)^j\langle x
\rangle^{-\sigma}_b \psi\| \leq O\left(\frac{1}{R}\right) \|g
x\left(\frac{x}{\langle x \rangle_b} \right)^j\langle x
\rangle^{-\sigma}_b \psi\|
$$

$$
\aligned & = O \left(\frac{1}{R}\right)\| g \left(\langle x
\rangle^{-1}_b x \right)^j
g^{-1} g x \langle x \rangle^{-\sigma}_b \psi \| \\
& \leq O \left(\frac{1}{R}\right)\| g\left(x\langle x
\rangle^{-1}_b\right)^j  g^{-1} \| \, \| gx \langle x
\rangle^{-\sigma}_b
\psi \|\\
&\leq O \left(\frac{1}{R}\right)\| \left( e^{i\beta} x\langle
e^{i\beta} x \rangle^{-1}_b\right) \|_{L^{\infty}_x}  \| gx \langle
x \rangle^{-\sigma}_b \psi \|.
\endaligned
$$
So, for $\beta$ sufficiently small $(R>1)$, the error terms from commuting $\langle
x \rangle^{-\sigma}_b $ are smaller than
$$
O \left(\frac{1}{R^2}\right) \|gx \langle x \rangle^{-\sigma}_b\psi
\|^2.
$$
Therefore, (4.16) implies
$$
\aligned
& p g^2(A/R) p + \langle x
\rangle^{-\sigma}_b x g^2(A/R) x \langle x
\rangle^{-\sigma}_b\geq \\
&\geq (1-\varepsilon) g(A/R) p^2 g(A/R) + (1-\varepsilon) g(A/R) \langle x
\rangle^{-\sigma}_b x^2\langle x
\rangle^{-\sigma}_b g (A/R)\\
&\geq \frac{1-\varepsilon}{4} g(A/R) \langle x \rangle^{-2\sigma}_b
g(A/R)
\endaligned
$$
where the last inequality follows from part b).
\enddemo
\demo{Proof of b)}
$$
\aligned
&p^2 + x^2 \langle x \rangle^{-2\sigma}_b =\\
&= \langle x\rangle^{-\sigma}_b p^2 \langle x\rangle^{-\sigma}_b +
\left(b^{-\sigma}- \langle x\rangle^{-\sigma}_b\right) p^2 \langle
x\rangle^{-\sigma}_b + \langle x\rangle^{-2 \sigma}_b x^2\\
&+ (b^{-\sigma}-\langle x\rangle^{- \sigma}) p^2(b^{-\sigma} - \langle
x\rangle^{- \sigma}) + \langle x\rangle^{-\sigma}_b p^2(b^{-\sigma}-\langle x\rangle^{-\sigma}_b)\\
&=\langle x\rangle^{-\sigma}_b
p^2 \langle x\rangle^{- \sigma}_b + (b^{-\sigma} - \langle
x\rangle^{- \sigma})p^2  (b^{-\sigma} - \langle
- x\rangle^{- \sigma}_b)\\
&+2 \sqrt{(b^{-\sigma} - \langle
 x\rangle^{- \sigma}_b)} \langle x\rangle^{-\sigma/2}_b p^2 \langle x\rangle^{-\sigma/2}_b\sqrt{(b^{-\sigma} - \langle
x\rangle^{- \sigma}_b)}\\
&+ 0 \left(\langle x\rangle^{-2\sigma -2}_b \left(\frac{x}{\langle
x\rangle_b}\right)^2\right) +x^2 \langle x\rangle^{- 2\sigma}_b
\endaligned
$$
which follows by commuting $\langle x\rangle^{-\sigma/2}_b$ and
$\sqrt{(b^{-\sigma} -\langle x\rangle^{-\sigma}_b)}$ through $p^2$.

For $b>> 1$, the results follows:
$$
\aligned & p^2 + x^2\langle x\rangle^{- 2\sigma}_b \geq \langle x
\rangle^{-\sigma}_b (x^2 + p^2)\langle x\rangle^{- \sigma}_b /2 +
x^2\left( \frac{1}{2}\langle x \rangle^{-2\sigma}_b - c \langle
x\rangle^{-\sigma -4}_b\right) \geq  \\
& \geq \frac{1}{2}\langle x\rangle^{- \sigma}_b(x^2 + p^2) \langle x\rangle^{-\sigma}_b \geq \frac{1}{4} \langle
x\rangle^{- 2\sigma}_b .
\endaligned
$$
\enddemo
\qed

\proclaim{ Theorem 4.4}

Let $V(x)$ be dilation analytic for all $|s|\leq \beta.$  Then
$$
i[V, \tanh A/R] = \frac{ + i}{2 \cosh A/R} \left\{ V^{[-\beta]} -
V^{[+\beta]}\right\} \frac{1}{\cosh A/R}
$$
where
$$
V^{[\beta]} \equiv e^{\beta A} V e^{-\beta A} =  V (e^{-i\beta} x).
$$
\endproclaim

\demo{Proof}
$$
\aligned & i\left[V, \frac{\sinh A/R}{\cosh
A/R}=\right]\frac{1}{\cosh A/R} i [V, \sinh A/R]
- \frac{1}{\cosh A/R}[V, \cosh A/R] \frac{\sinh A/R}{\cosh A/R}\\
& =\frac{1}{\cosh A/R}\left\{i[V, \sinh A/R] \cosh A/R - i[V, \cosh
A/R] \sinh A/R
\right\}\frac{1}{\cosh A/R}. \\
& \left\{\cdots\right\} = \frac{i}{4} \left[ [V, e^\beta] - [V,
e^{-\beta}] (e^\beta + e^{-\beta}) - \left( [V, e^{\beta}] + [V,
e^{-\beta}]\right)(e^\beta + e^{-\beta})\right]\\
& = \left[ 2[V, e^\beta] - 2[V e^{-\beta}] e^\beta\right]\frac{i}{4}\\
& = \left[ 2(V e^\beta- e^{+\beta} V) e^{-\beta} - 2  (V e^{-\beta} - e^{-\beta} V) e^{\beta} \right]\frac{i}{4}\\
& = \frac{i}{2}\left[V-e^\beta Ve^{-\beta} - V +  e^{-\beta} V e^\beta\right]\\
& = \frac{i}{2}\left[V^{[-\beta]} - V^{[\beta]}\right].
\endaligned
$$
\enddemo
\qed

\head Section 5\endhead
\subhead
5. Repulsive potentials and small Perturbations\endsubhead
\medskip

Let
$$H =-\Delta + V(x) + \varepsilon W (x)  \tag 5.1
$$
where $V, W$ as before, and have some analytic structure:

\noindent{\bf Assumption AN}

For some $\beta_0$ small, and $|\beta|\leq \beta_0$
$$
V(e^{\pm i \beta}x), W(e^{\pm i \beta} x)
$$
are bounded, continuously differentiable, and decay at $\infty$;
$$
V, x\cdot \nabla V, (x\cdot \nabla)^2 V, W, x\cdot \nabla W, (x\cdot \nabla)^2 W
$$
are all uniformly bounded by $C\langle x \rangle^{-2}$, and the
same holds for the analytic continuations(with $\vert \beta\vert \leq \beta_0)$ above.

\proclaim{Proposition 5.1}
Let $H$ as above, with $V, W$ satisfying Assumption AN.

Then
$$\aligned
&i[H, \tanh A/R] =\\
&= p g^2(A/R) p + \frac{1}{\cosh (A/R)} (i/2)
\left[ V(e^{i\beta} x) - V(e^{-i\beta} x)\right]\frac{1}{\cosh (A/R)}\\
&+ \varepsilon \frac{1}{\cosh (A/R)} (i/2)\left[ W(e^{i\beta} x ) -
W(e^{-i\beta} x )\right]\frac{1}{\cosh (A/R)}\\
& \beta = \frac{1}{R}.
\endaligned \tag 5.2
$$
\endproclaim

\noindent{\bf Definition}

$V$ is analytic repulsive potential if
$$
i \left( V(e^{i\beta} x) - V(e^{-i\beta} x )\right) \geq 0.
$$

\noindent{\bf Example}

$$
V(x) = \frac{1}{1+x^2} .
$$
In this case
$$
i \left(V (e^{i\beta} x) - V(e^{-i\beta} x)\right) = -2 \,  Im
\frac{1}{1+ e^{2i\beta} x^2}
$$
$$
= \frac{2 x^2 \sin 2\beta}{|1 + e^{2 i \beta} x^2 |^2} \geq c_\beta
\frac{x^2}{\langle x \rangle^4} , c_\beta > 0
$$
provided $|\beta| < \pi/ 4$.

We conclude that

\proclaim{ Theorem 5.2}

Let $H$ be as in (5.1), and $V,W$ satisfy the assumption AN.

Suppose, moreover that $V(x)$ is an analytic-repulsive potential, with lower bound
$$
i \left([V(e^{i\beta} x) - V(e^{-i\beta} x) \right] \geq c
x^2\langle x \rangle^{-\sigma}, c > 0,\quad \sigma \ge 4,
$$
and $W$ with decay of the above expression (to at least) of order
$\langle x \rangle^{-\sigma +2}$.  Then, for all $\varepsilon$ small
enough the RHS of equation 5.2 is positive and the corresponding
local propagation estimates hold:
$$
\int^T_0 \langle[ \| g(A/R) p \psi (t) \|^2
+\|\langle x \rangle^{-1} g (A/R) \psi (t) \|^2
+  \| x \langle x \rangle^{-2} g (A/R) \psi (t) \|^2\rangle] dt
$$
$$
\leq c|\langle \psi (T),   (\tanh A/R) \psi (T) \rangle |
+ |\langle \psi (0), \tanh (A/R) \psi (0) \rangle |.
$$
Here,
$$
g^2(A/R) \sim \frac{1}{R} \frac{1}{\cosh^2 A/R}.
$$
\endproclaim

\demo{Proof}
The only thing to check is that the $W$ term in the commutator, is bounded by the
repulsive contribution, coming from $-\Delta+V$.  To this end, note that near $x=0$,
$$
W(e^{i\beta} x) - W(e^{-i\beta} x) = x \int^\beta_{-\beta} e^{is}
W (e^{is} x) ds
$$
is $\sim x$.
\enddemo
\qed

\noindent{\bf Remark }

The condition of analyticity is technical, and is due to the fact that
the propagation observable we use is exponentially localized, up to a
constant, at $\infty$.

\head Section 6\endhead
\subhead
6. Local Decay and other propagation estimates \endsubhead
\medskip

The operator $\tanh A/R$ can play the role leading to an analytic
version of the projections on outgoing and incoming waves $P^\pm (A).$

We define
$$
F^+_M = F\left( \frac{A-M}{R}\right) = \left(\tanh \frac{A-M}{R} +
1\right) /2\tag 6.1
$$

So, $F^+_M$ is exponentially small (in $M/R$)  for $A- M < 0$.

Similarly, we define
$$
F_M^- = F^- \left( \frac{A+M}{R}\right) = \left(1-\tanh
\frac{A+ M}{R}\right) /2\tag 6.2
$$
We also notice the following inequality as a consequence of Thm 3.3,
Lemma 4.3.a, and proposition 5.1:

\proclaim{ Theorem 6.1}

For $H=-\Delta +V$ with $V$ satisfying assumption AN,  for all $R$
large enough, we have that:

$$
i[H, \tanh (A/R)] \geq (1 -\varepsilon) g_R(A) p^2 g_R(A) + g_R(A) \tilde
V_\beta g_R (A)
$$
where
$$
g^2_R(A) \sim \frac{2}{R} \frac{1}{2+ \cosh 2 A/R}\tag 6.3b
$$
and
$$
2\tilde V_\beta = i V(e^{i\beta} x) - i V(e^{-i\beta} x) \tag 6.3c
$$
$$
V(e^{i\beta} x) = e^{-\beta A} V (x) e^{+\beta A}  \tag 6.3d
$$
$$
\beta = 1/R.\tag 6.3e
$$
\endproclaim

It is now easy to find classes of potentials for which we get
monotonic decay estimates:

In one dimension we need either one of :
$$
\tilde V_\beta \geq 0, \text{ or } 2p^2 \sin 2\beta + \tilde V_\beta \geq 0,\tag i
$$
$$
\tilde V_\beta \geq  x^2 \langle x \rangle^{- 2\sigma}_b - \frac{1}{10} \langle x \rangle^{- 2\sigma+2}_b, \quad \sigma \ge2.\tag ii
$$
$$
V=V_1 + \varepsilon W  \tag iii
$$
where $V_1$ satisfies (ii) and $\varepsilon << 1$, and $|W_\beta| \leq
\langle x \rangle^{- 2\sigma+2}_b.$
$$
\align
&\text{Suppose that } -\Delta + V \geq 0.\tag iv \\
& \text{ Then, } p^2 + \tilde V_\beta/(2\sin 2\beta)
= ap^2 + (1-a) (p^2 +  V) + [\tilde V_\beta/(2\sin 2\beta) - (1-a) V]\\
& \text{ so, we need } \tilde V_\beta/(2\sin 2\beta) - (1-a) V\geq 0 \text{ for
some } 0 < a \leq 1.
\endalign
$$

\noindent{\bf In three dimensions }  Monotonic Decay estimates hold
whenever $p^2 + \tilde V_\beta/(2\sin2\beta) \geq 0$:

E.g., when,

$$
\frac{1}{4|x|^2} + \tilde V_\beta/(2\sin 2\beta) \geq 0, \text{ or when} \tag i
$$

$$
\aligned
& \text{Suppose that } H=- \Delta + V \geq 0. \\
&\text{ Then, we require that }\\
& p^2 + \tilde V_\beta/(2\sin 2\beta) = a p^2 + [(1-a) p^2 + \tilde V_\beta/(2\sin 2\beta)]\\
&= ap^2 + [-(1-a) V + \tilde V_\beta/(2\sin 2\beta)] + (1-a) (p^2 + V)\\
& \geq \frac{a}{4|x|^2} + [\tilde V_\beta/(2\sin 2\beta) - (1-a) V] \geq 0.
\endaligned
\tag ii
$$
which may be useful when $V$ has a negative part.
\medskip

\noindent{\bf Local Decay}

We have that for $F \equiv F(A/R) = \tanh A/R$
$$
i[H, F] = 2 pg^2(A) p + (1/\cosh (A/R ))\tilde V_\beta (1/\cosh (A/R ))\tag 6.5
$$
which we now assume to be positive: $\tilde V_\beta \geq 0$, and
$$
\aligned
& i[H, F] = 2 pg^2(A) p + (1/\cosh (A/R )) \tilde V_\beta (1/\cosh (A/R ))\\
& \geq 2 (1-\varepsilon) g(A) p^2 g(A) + (1/\cosh (A/R )) \tilde V_\beta (1/\cosh (A/R )) \geq \\
&\geq g(A) B^2 g(A) ,\quad \text{with}  B^2 \geq 0.
\endaligned
\tag 6.6
$$
Occasionally we have
$$
B^2 > \delta_{int} \chi (|x|\leq b_{int} ) + \delta_{out} |\tilde
V_\beta| \tag 6.7
$$
which is typical to one hump potentials $V$.

Now, let $M$ be a large positive number, and recall the definition:
$$
F^+_M (A/R) \equiv \left( \tanh\frac{A-M}{R} + 1\right)
/2
$$
and
$$
F^-_M(A/R) =\left(1-\tanh\frac{A+M}{R}\right)
/2
$$
the smooth projections on outgoing and incoming waves.

Then, letting for a momnet $f(A-M)\equiv 1/\cosh \frac{A-M}{R},$
$$
\aligned
& i[H, 2F^+_M] = 2p g^2(A-M) p+ f(A-M)\tilde V_\beta f(A-M)\\
& \geq 2 (1-\varepsilon) g(A-M) p^2 g(A-M) + f(A-M)\tilde V_\beta f(A-M)
\endaligned\tag 6.8
$$

$$
\aligned
& - i[H, 2F^-_M] = 2p g^2(A+M) p+ f(A+ M)\tilde V_\beta f(A+M)\\
& \geq 2 (1-\varepsilon) g(A+M) p^2 g(A+M)\\
& +f(A+M)\tilde V_\beta f(A+M).\endaligned
\tag 6.9
$$

In particular, it follows, since $-2 F^-_M \leq 0$ that

\proclaim{Proposition 6.2}
$$
\langle \psi(t), F^-_M \psi (t) \rangle \downarrow 0,\text{ as } t\to
+ \infty\tag 6.10a
$$
and
$$
\int^T_0 \|B g(A+M) \psi (t) \|^2 dt \leq \langle \psi(0), 2 F^-_M
\psi(0)\rangle, \tag 6.10b
$$
 with $B$ is defined in equation (6.6).
\endproclaim

This kind of monotonic decay is interesting, as it gives control of
the solution in the classically forbidden regions in terms of the size
of the solution  at {\bf time zero} with no corrections.

Applications will be discussed separately.

Next, we want to jack-up the decay estimate to a slowly decaying
weight, rather then $Bg$.

For this, we introduce new propagation observables:

$(\sigma >0)$
$$
\aligned
& 0\leq F^\pm_M(b^{-\sigma}-\langle x \rangle^{-\sigma}_b)+(b^{-\sigma}-\langle x \rangle^{-\sigma}_b)F^\pm_M \equiv F^\pm_M(b^{-\sigma}-\langle x \rangle^{-\sigma}_b) +c.c. \\
& \langle x \rangle^{-\sigma}_b
= (b^2 + |x|^2)^{-\sigma /2}\leq b^{-\sigma}.\endaligned
\tag 6.11
$$

We then have:(c.c. stands for Hermitian conjugate)

\proclaim{Proposition 6.3}
$$
\aligned
& i[H, F^+_M(b^{-\sigma} -\langle x \rangle^{-\sigma}_b) +c.c.]=\\
& =2\sigma \langle x \rangle^{-\frac{\sigma}{2}-1} A(F^+_M)\langle x
\rangle^{-\frac{\sigma}{2}-1}\\
& +\sum_i\tilde F_M C_i \tilde F_M + O (R^{-a}) F_M O(1) A \tilde
F_M
\endaligned
\tag 6.12
$$
where $\tilde F_M$ stands for approximate (discrete) derivatives of
$F_M$ (w.r.t. $A$), and $C_i, O(1),$ are operators which are of
higher order in $\langle x \rangle^{-1}$, and of order $R^{-1}$ at least,
$R$-large.
\endproclaim

\demo{Proof}
We denote $\langle x \rangle_b \equiv \langle x \rangle,$ and let $g(A)\equiv 1/\cosh\frac{A-M}{R}.$
$$
\aligned
& i[H, F^+_M(b^{-\sigma} -\langle x \rangle^{-\sigma}) +c.c.]=\\
&  =i[H, F^+_M] (b^{-\sigma} -\langle x \rangle^{-\sigma})  +
(b^{-\sigma} -\langle x \rangle^{-\sigma}) i[H, F^+_M]\\
& + F^+_M i [p^2,  -\langle x \rangle^{-\sigma}] +c.c.\\
& =g(A)(2\sin 2\beta p^2+\tilde V_\beta) g(A)(b^{-\sigma}-\langle x
\rangle^{-\sigma}) +  (b^{-\sigma}-\langle x
\rangle^{-\sigma}) g(A)(2\sin 2\beta p^2 +\tilde V_\beta )g(A)\\
& +F^+_M \sigma [ \langle x \rangle^{-\sigma-1}xp +  p x\langle x
\rangle^{-\sigma-2}] +c.c. \equiv I +I^* + J. \endaligned
\tag 6.13
$$
\enddemo
We symmetrize $J$ first:

Since
$$
\aligned
& A  = \frac{1}{2}(x \cdot p + p\cdot x) = x\cdot p - n i /2 = p\cdot x +
n i /2,\\
& J= \sigma F^+_M[\langle x \rangle^{-\sigma-2} A + A \langle x
\rangle^{-\sigma-2} ] +c.c.\\
& = \langle x \rangle^{-\sigma-2}\sigma A(F^+_M) + \sigma A(F^+_M)
\langle x \rangle^{-\sigma-2}\\
& =\sigma \left[ [AF^+_M, \langle x \rangle^{(-\sigma-2)/2}],\langle x \rangle^{(-\sigma-2)/2}
\right]+2\sigma \langle x \rangle^{(-\sigma-2)/2}AF_M^+\langle x \rangle^{(-\sigma-2)/2}.\endaligned
\tag 6.14
$$
We need to know that we can write
$$
[F(A), \langle x \rangle^{-2}] \text{ as } \tilde F(A) C
$$
with $C$ bounded, of order $\langle x \rangle^{-2}$, at least.

Now,
$$
[F(A), \langle x \rangle^{-2}] = - \langle x \rangle^{-2}[F(A),
x^2]\langle x \rangle^{-2} \tag 6.15
$$
$$
=-\langle x \rangle^{-2} 2x \tilde F(A) x \langle x \rangle^{-2}.
$$
Then, using that $\beta = \frac{1}{R}$ is small, we can write for any $\beta',$ $(g_{\beta{'}}(A)\equiv 1/\cosh(\beta{'} A)$
$$
\langle x \rangle^{-2} x \tilde F(A) x \langle x \rangle^{-2} = g_{\beta{'}}(A)
\cosh (\beta{'} A) x \langle x \rangle^{-2} \tilde F(A) x \langle x
\rangle^{-2} \cosh(\beta{'} A) g_{\beta{'}}(A)
$$
$$
\aligned
&=\frac{1}{2} g_{\beta^{'}} (A) (x \langle x \rangle^{-2})_{\beta^{'}}\cosh
2\beta^{'} A \tilde F(A) (x \langle x \rangle^{-2})_{\beta^{'}}
g_{\beta^{'}} (A)\\
&+\frac{1}{2} g_{\beta^{'}} (A) (x \langle x \rangle^{-2})_{-\beta^{'}}\cosh
(-2\beta^{'} A) \tilde F(A) (x \langle x \rangle^{-2})_{-\beta^{'}}
g_{\beta^{'}} (A)\\
&+\frac{1}{2} g_{\beta^{'}} (A) (x \langle x \rangle^{-2})_{-\beta^{'}}
 \tilde F(A) (x \langle x \rangle^{-2})_{+\beta^{'}}
g_{\beta^{'}} (A)\\
&+\frac{1}{2} g_{\beta^{'}} (A) (x \langle x \rangle^{-2})_{+\beta^{'}}
 \tilde F(A) (x \langle x \rangle^{-2})_{-\beta^{'}}
g_{\beta^{'}} (A)
\endaligned
$$
where
$$
( x \langle x \rangle^{-2})_\beta^{'} = e^{i\beta^{'}}
x\langle e^{i\beta^{'}}x\rangle^{-2}  = e^{i\beta^{'}} x ( 1+
e^{2i\beta^{'}} x^2)^{-1}
$$
$$
= e^{i\beta^{'}} x( 1+ e^{-2i\beta^{'}} x^2) (1+ e^{2i\beta^{'}}
x^2)^{-1} (1+e^{-2i\beta^{'}} x^2)^{-1}
$$
$$
2 Re \left[ e^{i\beta^{'}} x (1 + x^2\cos 2 \beta^{'} - x^2 2 i \sin
2 \beta^{'}) ( 1+x^4 + 2 x^2 \cos 2 \beta^{'})^{-1}\right] = O
(\langle x \rangle^{-1})
$$
and similarly for the Imaginary part,
(for $\beta^{'}$ small).  Here we choose $\beta^{'}\leq \beta$.
$$
\tilde F(A) \sim O \left(\frac{1}{R}\right)(\cosh 2\beta A)^{-1},
$$
$$
g_{\beta^{'}}(A) = (\cosh \beta^{'} A)^{-1}.
$$
So, we have that
$$
[F(A) , \langle x \rangle^{-2}] = g(A) Cg(A).\tag 6.16
$$
Similarly, we can rewrite
$$
[F(A) , \langle x \rangle^{-2}] =-2\langle x \rangle^{-2} \tilde
F_-(A) x^2 \langle x \rangle^{-2}
$$
$$
= +2 \langle x \rangle^{-2}\tilde F_-(A) \left( \frac{1}{1+ x^2}
-1\right)
$$
$$
=-2 \langle x \rangle^{-2} \tilde F_- (A) + 2 \langle x
\rangle^{-2}\tilde F_-(A) \left(\frac{1}{1+ x^2}\right)(\cosh
\beta^{'} A)(\cosh \beta^{'} A)^{-1}
$$
$$
=-2 \langle x \rangle^{-2} \tilde F_- (A) + 2 \langle x
\rangle^{-2}\left\{ \tilde F_- (A) e^{\beta^{'} A}\left(\frac{1}{1+
x^2}\right)_{\beta^{'}} + \tilde F_-(A) e^{-\beta^{'} A}
\left(\frac{1}{1+ x^2}\right)_{\beta^{'}}\right\} (\cosh \beta^{'}
A)^{-1}
$$
$$
=\sum C_i \tilde F(A) , \qquad \tilde F(A) \sim (\cosh \beta^{'}
A)^{-1} \tag 6.17
$$
$$
C_i = O (\langle x \rangle^{-2} \frac{1}{R}) ,\qquad  \beta^{'} \leq
\beta, \text{ small}.
$$

Using the above identities for $[F(A) , \langle x
\rangle^{-2} ]$ we can easily symmetrize the expressions for $I, I^*$
and $J$ to get:

$$
J= \langle x \rangle^{-\frac{\sigma}{2}-1} 2 \sigma A(F^+_M)
\langle x \rangle^{-\frac{\sigma}{2}-1}
$$
$$
+\left[\langle x \rangle^{-\frac{\sigma}{2}-1}, [\langle x
\rangle^{-\frac{\sigma}{2}-1}, \sigma A(F^+_M)]\right]
$$
$$
= \langle x \rangle^{-\sigma/2-1} 2\sigma A(F^+_M) \langle x
\rangle^{-\frac{\sigma}{2}-1}\tag 6.18
$$
$$
+\left[\langle x \rangle^{-\frac{\sigma}{2}-1}, [\langle x
\rangle^{-\frac{\sigma}{2}-1}, \sigma A(F^+_M)]\right]
$$
Using that for any $Q,$
$$
\left[[Q, f(A)], g(A)\right] = [Q,g] f- f[Q,g] = \left[[Q,g],
f\right]\tag 6.19a
$$
$$
i[\langle x \rangle^{-m}, AF^+_M] = +m\langle x \rangle^{-m-2} x^2 F^+_M
+ A i [\langle x \rangle^{-m}, F^+_M]\tag 6.19b
$$
$$
= m\langle x \rangle^{-m-2} x^2 F^+_M - AC\tilde F_M = m \langle x
\rangle^{-m-2} x^2 F^+_M - C A \tilde F_M - [A, C]\tilde F_M
$$
$$
[A, C] = O (\langle x \rangle^{-m}/R).
$$

Commuting again with $(m\equiv \sigma/2+1)$ $\langle x
\rangle^{-\frac{\sigma}{2}-1}$ we get that the double commutator is of
the form:
$$
  O (\langle x \rangle^{(-\sigma-2)/2}/R)
\tilde F_M .O (\langle x \rangle^{(-\sigma-2)/2}/R).
$$
Therefore
$$
J= \langle x \rangle^{-\sigma/2-1} 2\sigma A(F^+_M) \langle x \rangle^{-\frac{\sigma}{2}-1}
$$
$$
\aligned
&+ O(R^{-1}) (\langle x \rangle^{-\frac{\sigma}{2}-1} \tilde F_M)
O(1) F_M \langle x \rangle^{-\frac{\sigma}{2}-1}\\
&\ge\sigma\ \langle x \rangle^{-\frac{\sigma}{2}-1}\vert A \vert F_M^+\langle x \rangle^{-\frac{\sigma}{2}-1}\endaligned
\tag 6.20
$$
Symmetrizing $I + I^*$, we have that, as above:

$$
\aligned & I + I^* =\tilde F_M H_{\beta} \tilde F_M\chi_b^2(|x|)+\chi_b^2(|x|)\tilde F_M H_{\beta} \tilde F_M\\
&= \tilde F_M 2 \chi_b(|x|) (p^2 + \tilde V_\beta)
\chi_b(|x|) \tilde F_M + \tilde F_M O (\langle x \rangle^{-\sigma-2}
) \tilde F_M
+\\
&  \tilde F_M O(\langle x \rangle^{-\sigma} R^{-1})H_{\beta}\tilde F_M+c.c.\\
&\ge \tilde F_M \chi_b(|x|) (2\sin 2\beta p^2 + \tilde V_\beta)
\chi_b(|x|) \tilde F_M ;
\endaligned\tag 6.21
$$

$$
\chi_b(|x|) = (b^{-\sigma} -  \langle x \rangle^{-\sigma}_b)^{1/2}
$$
Combining (6.20), (6.21) we have that:

\proclaim{ Theorem 6.4}  (Local Decay for Analytic Repulsive
Potentials)

Let $H =-\Delta + V(x)$ as before and s.t. $V$ is Analytic
repulsive, and $-\Delta +V \geq 0$.

Then, for $\sigma >0,$
$$
\int^T_0  \|\langle A\rangle^{1/2} F^+_M \langle x
\rangle^{-\sigma - 1} \psi \|^2 dt
$$
$$
+\int^T_0  \| p\chi_b (|x|) \tilde F_M \psi \|^2 dt \leq C\|\psi\|^2
$$
\endproclaim

\noindent{\bf Remark:}

We can replace $AF^+_M $ by $\langle A\rangle$ in the expression for
$J$, eq (6.20), using the local decay estimate proposition (5.2),
which controls the region $|A| \leq $C, and a  similar bound on $F_M^-.$

Similar estimate holds for $F^-_M$:

\proclaim{ Theorem 6.5}  (Pointwise (and integral) decay of Incoming
waves)

Under the conditions of Theorem 6.4, we have that
$$
\int^T_0\|\langle A\rangle^{1/2} F^-_M \langle x
\rangle^{-\sigma -1} \psi \|^2 dt + \langle \langle x
\rangle^{-\sigma}\psi (T),
(F^-_M)^ A \langle x
\rangle^{-\sigma}\psi (T)\rangle
\leq 2 \langle \psi (0) (F^-_M)^2 \psi(0)\rangle .
$$
\endproclaim

Combining all the above, we get that local decay holds with the
following weight:
$$
\int^T_0 \|\langle A\rangle^{1/2}\langle x
\rangle^{-1-\varepsilon} \psi \|^2 dt \leq c\|\psi\|^2.
$$

\head Section 7\endhead
\subhead
7. Applications:  Schwarzschild manifolds, generalized repulsive
   potentials\endsubhead
\medskip

When the Hamiltonian $H\geq 0$, we can get the desired decay estimates
by simply verifying that
$$
 2\sin 2\beta p^2 + V_\beta > 0
$$
for some $\beta$ small.

In particular, if $- x\cdot\nabla V > 0$, together with some
uniformity of the analytic continuations $V_\beta,$ the above
inequality follows.

We also get local decay, for one hump potentials, including the
Schwarzschild for each fixed angular momentum:
\medskip

Case Study: Schwarzschild potentials

Here we solve the wave equation
$$
-\frac{\partial ^2 u}{\partial t^2} = Hu
$$
$$
u_0 = (f_0,g_0) \in H^1 \otimes L^2
$$

Let
$$
H = - \partial^2_{r_*} + V_\ell (r) \text{ on } L^2 (\Bbb R, dr_*)\tag
7.1
$$
where $r_*= r + 2 M \ln(r - 2M)$

so that $\frac{dr_*}{dr} = 1 + 2M \frac{1}{r-2M} = \frac{r- 2 M +
2M}{r-2M}$
$$
= \frac{r}{r-2M}\text{ and } \frac{dr}{dr_*} = \frac{r-2M}{r} = 1 -
\frac{2M}{r}.\tag 7.2
$$
$$
V_\ell (r) = (1-\frac{2M}{r} )\frac{2M}{r^3} + \left( 1-\frac{2M}{r}\right)
\frac{1}{r^2} \ell (\ell +1) \quad \ell = 0, 1,2,\ldots .\tag 7.3
$$
Since for each $\ell, V_\ell(r)$ is a one hump potential around the
point
$$
\left(x^2 = \ell(\ell + 1)\right) \,  r_\ell  \equiv
\frac{3M(\lambda^2 -1) + M\sqrt{\rho(\lambda^2 -1)^2 + 3 2 \lambda^2}}
{2 \lambda^2}
$$
$$
\alpha^*_\ell = \alpha^*(r=r_\ell) ; \alpha^*_\infty = r_*(r= 3M) \tag
7.4
$$
it follows that the decay estimates hold for analytic each $H_\ell$,
if we can show that the humps are repulsive!  Summing over all
$\ell$, after multiplying by $P_\ell$, the projection on the
$\ell$'th spherical harmonic, local decay follows for $-\Delta$ on
Schwarzschild manifolds.

This argument applies to all manifolds where the resulting potential
is one-hump, analytic repulsive at fixed angular momentum.

In fact we get somewhat different and new estimates in this case,
since, as we remarked before, the propagation observable(PROB) we use is
bounded on $L^2$, unlike the Morawetz estimate and its various
generalizations which are bounded from $H^{1/2} \to L^2$.
 The solution of the wave equation can be written in terms of the
 initial data $ u(x,t=0):=f_0, \dot u(x,t=0):=g_0$ as:
 $$
 u(x,t):={\text U}(t)u_{0}=\cos (\sqrt H t)f_0+\frac{\sin \sqrt H t}{\sqrt
 H}g_0.
 $$
There is a fundamental new difficulty with the WE (Wave Equation) as compared with
the Schr\"odinger equation. This is due to the fact that $L^2$ norm
can grow linearly in time for the WE, and the LHS of the propagation estimate(PRES) has a form
different from the Schr\"odinger case.

\proclaim{Theorem 7.1}{Local Decay-WE}
$$
\int_0^T \| \langle x \rangle ^{-3/2} \langle A \rangle ^{1/2}u \|^2
dt < C E^{1/2}(u_0)\left(E^{1/2}\left( g (H\leq \epsilon) u_0\right)+c_{\epsilon}\|u_0\|_{L^2}\right ). \tag
7.5
$$
\endproclaim

The proof of the Theorem is a consequence of the propositions
that follow:

The Heisenberg equation formulation of the wave equation is
$$
\partial_t[(u, B\dot u) - (\dot u, Bu)] = [H, B]\tag 7.6a
$$
where
$$
-\partial_t^2 u = Hu.\tag 7.6b
$$
Using $B\equiv \frac{\partial}{\partial t},$
we get the {\bf Energy Identity}:
$$
\aligned
&\frac{\partial}{\partial t} [(u,\frac{\partial}{\partial t}\dot u)-(\dot u,\frac{\partial}{\partial t}u)]\\&=
\frac{\partial}{\partial t}[(u,-Hu)-(\dot u,\dot u)]\\&=
\frac{\partial}{\partial t}\int |\nabla u|^2 +|\dot u|^2 +V(x)|u|^2 dx \equiv \frac{\partial}{\partial t}E(u)=0.
\endaligned
$$
So, the energy conservations reads $E(u)=E(u_0).$
In our case
$$
B=i \tanh A/R.\tag 7.6c
$$
First, we reduce the problem to initial data with localized
frequencies near zero.

For this, let $g = g(|p|\leq 1), \bar g = 1 -g$ and write $u$ as
$$
\aligned
& u = gu + \bar g u\\
& (u, F\dot u) - (\dot u, Fu) = (gu, F g\dot u) + (\bar g u, Fg\dot u)
\\
&+(\bar gu, F\bar g\dot u) + (gu, F\bar g\dot u) - (g\dot u, Fgu)
-(\bar g\dot u, F gu)\\
&-(\bar g \dot u, F\bar g u) - ( g\dot u, F\bar g u).
\endaligned
$$
Every term with $\bar gu$ is good.

$$
\bar g u = \bar g |p|^{-1} \langle p\rangle \langle p\rangle^{-1} p u
$$
and therefore $|(\psi, \bar g u)| \leq \|\psi\| \,  \|\bar g
|p|^{-1}\langle p \rangle \| \,  \|\langle p\rangle^{-1} pu\|$.

Next, we have
$$
\aligned
& (g u, F\bar g \dot u) - ( \bar g\dot u, Fg u) = \\
&=(pgu, F_({\frac{1}{p}})\bar g \dot u) + ( p gu, \tilde F_{+} \frac{1}{p}
\bar g \dot u)\\
& -(\frac{1}{p}\bar g \dot u, Fp g u) - (\frac{1}{p} \bar g\dot u,
 \tilde F_- pgu)\\
&\leq 2 \|\frac{1}{p} \bar g \| \, \|\dot u\|(\|F\|+ 2 \|\tilde F\|) \|g
 pu\|.\endaligned
$$
Finally to deal with terms with no $\bar g$ in them, we need to
exploit the fact that  $\tanh \, A/R$ vanishes linearly in $A$ near zero.

$$
\aligned
& -(g \dot u, F gu) + ( gu, F g\dot u) = \\
&-\left(  g\dot u, F \frac{1}{A} (x p- i/2) g u\right) + \left( ( x p -i/2) gu,
\frac{1}{A} Fg\dot u\right)\\
& =-\left(g\dot u, (F\frac{1}{A} x \tilde g) gpu\right) + (g pu, \tilde g x
\frac{1}{A} F g\dot u) ,\qquad \tilde g g = g.\endaligned
$$

Now, since $F= i \tanh A/R,\qquad \pm i F =\mp G$ with $G=G^*$.  Furthermore,
$g$ acts like the convolution with the function $\hat g$, the Fourier
transform of $g$, which is real.

Hence $gu, g\dot u$ are real.

This leads to the cancellation of the two terms with $-\frac{i}{2}$
factor.

We are therefore left with
$$
\big|2(g\dot u, (FA^{-1} x \tilde g) g pu) + 2 (g pu, (\tilde g
 x\frac{1}{A} F) g\dot u)\big|
$$
$$
\leq 4\|g\dot u \|\|FA^{-1}x\tilde g gpu\|.
$$

Hence, collecting all the terms, we arrive at
$$
\big|\langle u, i \tanh A/R \dot u\rangle - \langle \dot u, i \tanh A/R
u\rangle \big|
$$
$$
\leq c\|\dot u \| \, \| \langle p \rangle^{-1} pu \|+c\|g\dot
u\|\|FA^{-1}xgpu\|.\tag 7.7
$$

To this end, we need the following propagation observables, and energy decomposition;
Fix a (large) time T.
We break the initial data  $(f,g) = u_0$ as:
$$
\aligned
 u_0 &= F(H\leq T^{-1}) u_0 + F(H\ge T^{-1})u_0 \\
 &= u_l+u_h \equiv F_{<} u_0 +F_{>} u_0
 \endaligned
$$
Clearly then, since H commutes with the dynamics $ {\text U}(t)$,
$$
$$
 that
$$
 {\text U}(t)u_{_l} = F(H\leq T^{-1}) {\text U}(t)u_{_l}
$$
$$ {\text U}(t)u_{_h} = F(H\geq T^{-1}) {\text U}(t) u_{_h},
$$
so that:
$$
\aligned
\|H{\text U}(t)u_{_l}\|_{L^2}&=\|HF_<{\text U}(t)u_{_l}\| \\
& \leq T^{-\frac{1}{2}}\|H^{\frac{1}{2}}{\text U}(t)u_{_l}\|_{L^2}
E^{\frac {1}{2}}(u_{_l}).
\endaligned
$$
We will use the following propagation observables:
$$
\aligned
& B_1 = i\tanh( A/R) \\
& B_2 = iF^{\pm}_M(A/R) \\
& B^\sigma_n\equiv F_M^\pm i \langle x \rangle ^{-m} +c.c.
\qquad  m \geq 0.
\endaligned \tag 7.8
$$
We then have, as before, that
$$
\aligned
& [H,B_1] = \tilde g_{_0}(A) H_{_{\beta}}  \tilde g_{_0} (A) \\
& [H,B_2] = \pm \tilde g _{M} ^{\pm} (A) H_{_{\beta}} \tilde g _M
^{\pm}(A) \\
& [H,B_m ^{\sigma}] = F_{M} ^{\pm}  \left \{ \langle x \rangle
^{-m-2} A + c.c. \right \} F_M ^{\pm} \pm \left \{ \tilde g _M
^{\pm} H_{_{\beta}} \tilde g _{_M} ^{\pm } \langle x \rangle ^{-m}
F_M ^{\pm} + c.c. \right \}
\endaligned \tag 7.9
$$
Next, we have the following preliminary estimates on the $LHS$ of
the Heisenberg identity:
$$
$$
\proclaim{Lemma 7.2}
$$
\tag 7.10
$$
$$
 |\langle \dot u_h,B_1 u_h \rangle - \langle u_h, B_1 \dot u_h
\rangle | \leq C \| \dot u_h \| _{L^2} E^{\frac{1}{2}}(u_h)
T^{\frac{1}{2}} \tag i
$$
$$
 |\langle \dot u,B_1 u \rangle - \langle u, B_1 \dot u
\rangle | \leq C\| \chi(|A| \leq \lambda ) \dot u \| _ {L^2}
 \| \chi (|A| \leq \lambda ) (\tanh A/R) u \|_{L^2} \tag ii
$$
$$
\qquad \quad + C \| \chi ( |A|\geq \lambda) e^{-|A/R|} \dot u \|
_{L^2} \| \chi(|A| \geq \lambda) e^{-|A/R|} u \|_{L^2}
$$
$$
\text {Similar bounds hold for } B_1 \text {,with } A \to A-M. \tag
iii
$$
$$
|\langle  \dot u, B_m ^\sigma u \rangle - \langle u, B_m ^\sigma
\dot u \rangle | \leq CE(u_0)  \qquad \text {for} \qquad m \geq 1.
\tag iv
$$

\endproclaim
$$
$$
\demo {Proof}
$$
\text {Follows by Cauchy-Schwarz inequality and}   \tag i
$$
$$
 \|u_{_h}\|_{L^2} = \|H^{-\frac{1}{2}} F(H\geq \frac{1}{T}) H^{\frac{1}{2}} u_{_h} \| \leq T^{\frac{1}{2}} E^{\frac{1}{2}} (u_{_h}).
$$
$$
\text {Follows by noting that on support} \chi_{_{\gtrless}}(A):
\tag ii
$$
$$
|(\tanh (A/R) \mp 1 ) \chi_{_{\gtrless}} (A)|\lesssim 2e^{(-2|A/R|)}
\chi_{_{\gtrless}}(A)
$$
and that
$$\langle \dot u, Bu \rangle - \langle u, B\dot u \rangle
\equiv \langle B \rangle _u ^ {Heis} = 0
$$
$$\text {for } B=1  \text{  ( or any reality preserving symmetric
operator).}
$$
$$
\text{Follows from (i), (ii)  by replacing } A \text {  by  } A- R.
\tag iii
$$
$$
\text { Follows from } \| \langle x \rangle ^{-1} u \| _{L^2} \leq
CE^{\frac{1}{2}} (u), \text { and that} \tag iv
$$
$$
 [\langle x\rangle ^{-1}, F_{M} ^{\pm}]= -\langle x\rangle ^{-1} [\langle x \rangle, F_{M} ^{\pm}] \langle x \rangle ^{-1}
\simeq \langle x \rangle^{-1} \langle x \rangle \tilde F_{M} ^{\pm}
\langle x \rangle ^{-1} ,
$$
$$
\text { with } \tilde F_M ^{\pm}  \text { -  bounded}.
$$
\enddemo
We proceed to estimating $\langle B_1 \rangle _{u_{_h}} ^{Heis} $.

\proclaim{Proposition 7.3}
$$
$$
There exists a sequence of times, $ T_n \to \infty $ , such that
$$
\|F(|x| \leq MT_n^{\frac{1}{2}}) \tilde g(A) u_{_h}(T_n) \|_{L^2}
\leq CT_n^{\frac{1}{4}} M.\tag 7.11
$$
\endproclaim
\demo{Proof}
$$
\text { Applying the previous propagation estimates with } B_1,
\text {and using Lemma (7.2)} (i), \text {it follows that}:
$$
$$
\int_0^T \|\tilde g (A) \langle x \rangle ^{-1} u_{_h} \| ^2 dt +
\int _0 ^T \| \tilde g (A) pu_{_h}(t)\|^2 dt \leq CT^{\frac{1}{2}}
E(u_{_h}).
$$
Next, we apply the cutoff in  $|x|:$
$$
\|F(|x| \leq M\sqrt T) \tilde g (A) u_{_h}(t) \| \leq M \sqrt T \|
\langle x \rangle ^{-1} \tilde g (A) u_{_h}(t) \|,
$$
so that
$$
\aligned \int_0^T \| F(|x| &\leq MT^{\frac{1}{2}}) \tilde g(A)
u_{_h}(t) \| ^2 dt \\
& \leq M^2 T \int_0 ^T \| \langle x \rangle ^{-1} \tilde g(A) u_{_h}
(t) \| ^2 dt \leq CM^2T^{\frac{3}{2}}.
\endaligned
$$
Therefore, $\exists  T_n \to \infty \quad s.t. $
$$
\| F(|x| \leq M \sqrt {T_n}) \tilde g (A) u_{_h} (T_n ) \|^2  \leq
CM^2T_n^{\frac{3}{2}} T_n^{-1},
$$
$$
\text { since, otherwise}
$$
$$
\int_0^{T_n} \| F(|x| \leq M \sqrt T) \tilde g(A) u_h(t) \|^2 dt >
CM^2 T_n ^{\frac{3}{2}}.
$$

 Next, we use the above proposition to bound  $ \langle B_1 \rangle _{u_{_h}} ^{Heis}$,  using the fact that
 if $|x| > M\sqrt T $, and $ |p| \geq \frac{1}{\sqrt{T}}$,  then, classically, in the phase space, $ A\gtrsim M, $
 which, together with the localization in $ A$, via $ \tilde g(A)$, gives fast decay in $ M$, for $ |x| > M \sqrt T$.
\enddemo

\proclaim{Proposition 7.4}
$$
|\langle B_1 \rangle _{u_{_h}} ^{Heis} | \leq C \| F(|x| \leq
MT^{\frac{1}{2}}) \tilde g(A) u_{_h} (t) \| \|\dot u \| _{L^2}
$$
$$
+ O(M^{-\infty}) T E(u).
$$
\endproclaim
\demo{Proof}
$$
$$
We need to bound
$$
\aligned
\| F(|x| &> MT^{\frac{1}{2}}) \tilde g ^2(A) F(H\geq T^{-1}) u_{_h}(t)\|_{L^2}\\
&\leq TE^{\frac{1}{2}}(u_{_h})\|F(|x|>MT^{\frac{1}{2}})\tilde g^2(A)
F(H\geq T^{-1})\|_{L^2 \to L^2}
\endaligned
$$
where we used that
$$
\|u_{_h}(t) \| \leq CtE^{\frac{1}{2}}(u_{_h}).
$$
To this end, we write the above operator product as
$$
\|F( |x|>MT^{\frac{1}{2}}) x^{-2n}x^{2n}\tilde g^2(A) H^n H^{-n}
F(H\geq \frac{1}{T}) \|
$$
$$
\leq \|  x^{2n}H^n\tilde g^2(A) \| M^{-2n}T^{-n}
$$
$$
+\|x^{2n}[\tilde g^2(A), H^n]\|M^{-2n}
$$
The first term on the $RHS$ has a factor $x^{2n}H^n$, which , when
expanded, is a sum of terms of the form $x^{2n}P^{2j}V^k$ ... and
such that the order in $x$ is at most $2n-2j$ and the order in $p$
is $2n-2j$ in each monomial $P_j\thicksim x^{2n}V^kp^j...$.$$$$
 This is because
our $V(x)$ decays at least like $|x|^{-2}$.
 Hence, we can always
pair each monomial to be
$$
P_j \sim x^{2n-2j}p^{2n-2j} \left ( 1+O( {\frac{1}{x}})\right ) C_j.
$$
Each such monomial can be rewritten as
$$
P_j \sim C_j ^{'}\left (1+O(\frac{1}{|x|}) + O(A^{-1})\right )
A^{2n-2j}
$$
Hence the first term on the $RHS$ is bounded by
$$
\sum_j C_j^{''}\|A^{2n-2j} \tilde g_{_1}(A) \|.
$$
The second term on the $RHS$ is similar:
$$
x^{2n}[\tilde g^2(A), H^n] \sim \sum_j x^{2n} [\tilde g ^2 (A), P_j]
\sim \sum_k g_{k^{'}}(A)x^{2k}p^{2k} \sim \sum _l g_{_l}(l) A^{2l}.
$$
Using the exponential bound on the $\tilde g^2(A)$, and noting that the number of terms is at most of order $ n^n ,$
we get a bound of the form ( after inserting $x^NH^{N/2}$)
$$
M^{-N}\sum_{n=1} ^N n^n A^n e^{-|A/R|} \lesssim \left (
\frac{N^2R}{2eM}\right )^N.
$$
If we choose $N^2 \sim R^{-1} \sqrt M,$ we get a bound $\sim
M^{-N/2} \sim M^{-\frac{1}{2R}M^{\frac{1}{4}}}.$
\qed \enddemo
The propositions above imply ( after choosing $M \gtrsim ClnT)$.
\proclaim{Theorem 7.5}
$$
\left | \langle B_1 \rangle _{u_{_h}} ^{Heis} (T_n) \right | \leq
CT_n ^{\frac{1}{4}} ln T_n E(u_{_h}).
$$
\endproclaim
The above process, beginning with the bound of Lemma (7.2)(i), is
now iterated ( $ ln T$ times...), where we use the above
$T^{\frac{1}{4}}$ bound to replace the $T^{\frac{1}{2}}$ bound of
Lemma(i). This will give a $T^{\frac{1}{8}}$ bound etc... $$$$ We
conclude that

\proclaim{Theorem(7.6)}
$$
\left \langle B_1(T_n) \right \rangle_{u_{_h}}^{Heis} \leq
CE(u_{_h}), \qquad T_n \to \infty.
$$
\endproclaim
\smallskip

Next, we need to bound $\langle B_1 \rangle _{u_l} ^{Heis}$.
$$
$$
The method is similar to the previous case; however, the propagation
observables used need to be iterated, and the argument is a bit more
involved.
$$
$$
To this end we consider the part of the data where $ H\leq T^{-1} ,$ that
is, estimating $ u_{_l} $.

\smallskip

 In this case the propagation observable we use is of the
general form
$$
F_M^{\pm} i\langle x\rangle ^{-\sigma} F_M ^{\pm} = B_{\sigma , M}
^{\pm}. \qquad \sigma \geq 1.
$$
The commutator with $H$ has two parts; one comes from $\langle x
\rangle ^{-\sigma} $ and another from $F_{M} ^{\pm}$. They have
apposite sign, and therefore, we need to control one of them in
terms of the other.
$$
$$
Since $\sigma \geq 1 $, the $LHS$ of the Heisenberg equation is
uniformly bounded in time:
$$
\langle B_{\sigma,M}^{\pm} \rangle _{u_{l}}^ {Heis}\leq CE(u_{l}).
$$
We have that
$$
\aligned
[H,B_{\sigma, M} ^{\pm}] &= -\sigma F_M ^{\pm} \langle x \rangle ^{-\sigma -2} A F_M ^{\pm} + c.c.\\
& \pm [g_M^{\pm} H_{_\beta} g_M ^{\pm} \langle x \rangle ^{-\sigma}
F_M ^{\pm} + c.c.] \equiv C_\sigma + D_\sigma.
\endaligned
$$
Our goal is to show that in some sense $D_\sigma$ is higher order,
so the $C_\sigma $ term will give a propagation estimate. We iterate
on $\sigma $ to get the final bound. $$ $$
Symmetrizing $C_\sigma$ ,
as before, we get
$$
\aligned
C_\sigma &= -\sigma  \langle x \rangle ^{(-\sigma-2)/2} A F _M ^{\pm}
 \langle x \rangle ^{(-\sigma -2)/2}  + O ( \langle x \rangle ^{-\sigma -2} \tilde g (A)/R) \\
&= -\sigma \langle x \rangle ^{(-\sigma-2)/2} A F_M^{\pm} \langle x \rangle ^{(-\sigma-2)/2} +
{\text O} (\langle x \rangle ^{-\sigma-2}\tilde g(A)/R) \\
\endaligned
$$
$$
\aligned \left | \langle u_{_l},D_\sigma u_{_l} \rangle \right |
&\leq C \| \tilde g_{_M}(A) \langle x \rangle ^{-\sigma} u_l
\|\text{ }
\|H_{_\beta} \tilde g_{_M}(A) F(H\leq 1/T) u_{_l} \| \\
& \leq C\langle\{ \| \tilde g_{_M}  H_{\beta} F_< u_{_l}\|  +
\|[p^2,\tilde g_{_M}]F_<u_{_l}\|+
 \|V_\beta \langle x \rangle ^2 \| \text{ }   \|[\langle x \rangle ^{-2}, \tilde g_{_M}]
 F_{_<}u_{_l}\|\rangle\}\|\tilde g_{M}\langle x\rangle ^{-\sigma}u_{_l}\| \\
 & \leq C\left [\frac {T^{-1/2}}{\sqrt R} \| H_{_\beta} ^{1/2} u_{_l}\| + \|g_{_{1,M}}(A) p^2F_{_<}u_{_l}\|\right ]\|\tilde g_{M}\langle x\rangle ^{-\sigma}u_{_l}\|\\
& \qquad \qquad \qquad \qquad  + C\frac{1}{\sqrt R}\| \langle x
\rangle ^{-2} F_{_<}u_{_l}\|\,
\|\tilde g_{M}\langle x\rangle ^{-\sigma}u_{_l} \| \\
& \leq \frac{C}{\sqrt R}T^{-1/2} E^{1/2}(u_{_l})\|\tilde
g_{M}\langle x\rangle ^{-\sigma}u_{_l}\|+\frac{C}{\sqrt R}\|\tilde g_{M}\langle x\rangle ^{-\sigma}u_{_l}\| O(L^2(dt)).
\endaligned
$$
using that
$$
\langle x \rangle ^{-2} \leq CH_\beta
$$
$$
p^2 \leq CH_\beta.
$$
$$
\aligned \|\tilde g_{_M}(A) \langle x \rangle ^{-\sigma} u_{_l}\| ^2
&= \langle \tilde g_{_M}(A) \langle x \rangle ^{-\sigma} u_{_l},
\tilde g
_{_M}(A) \langle x \rangle ^{-\sigma} u_{_l} \rangle \\
& \leq \frac{C}{R} | \langle u_{_l},C_\sigma u_{_l}\rangle | +
Ce^{-M/2R} \| \tilde g(A)\langle x \rangle ^{-\sigma} u_{_l} \| ^2
\endaligned
$$
Provided
$$
\sigma \geq (\sigma + 2)/2
$$
$$
\tilde g(A) \sim (\cosh (A/R)) ^{-1}.
$$
Putting it all together, we have:
\proclaim{Proposition 7.7}

For $\sigma \geq (\sigma + 2)/2,$
$$
\aligned
|\langle u_{_l}, [H,B_{\sigma , M} ^\pm]u_{_l} \rangle | &\geq |\langle u_{_l}, C_\sigma u_{_l} \rangle | \\
& -\frac{C}{M^{1/2}R}T^{-1/2}E^{1/2}(u_{_l})| \langle
u_{_l},C_\sigma u_{_l} \rangle | ^{1/2}
\\
&- Ce^{-M/2R} \| \tilde g (A) \langle x \rangle ^{-\sigma} u_l \|
T^{-1/2}E^{1/2}(u_{_l})
\endaligned
$$
For $B_{\sigma,M} ^{-} , \langle u_l,C_\sigma u_{_l} \rangle $ is
positive, and $\langle u_{_l}, C_\sigma u_l \rangle $ is negative
for $B_{\sigma,M} ^{+}$.
\endproclaim

We integrate over time the Heisenberg
equation and using the above proposition to obtain the following
propagation estimate:

\proclaim{Proposition 7.8}
$$
$$
For $ \sigma \geq \frac{\sigma}{2} + 1$ :
$$
\aligned
\int_0^T\|\langle A \rangle ^{1/2} \tilde F_M^{\pm} \langle & x \rangle ^{(-\sigma-2)/2} u_{_l}\|^2 dt \\
& \leq CT^{-1/2} E^{1/2}(u_{_l}) \int_0^T \| \tilde F_M^{\pm} \langle x \rangle ^{(-\sigma -2)/2} u_{_l}\|ds \\
& + CT^{-1/2}e^{-M/2R} \int_0^T\| \tilde g(A) \langle x \rangle
^{-\sigma} u_{_l}\|ds + CE(u_{_l})
\endaligned
$$
where the last term comes from  $\langle B_\sigma \rangle _{u_{_l}}
^{Heis} $.
\endproclaim
$$
$$
Applying the above result with $\sigma = 2$, we can
get the following local decay estimate:

\proclaim{Theorem7.9}
$$
$$
Under the previous assumptions on the
Hamiltonian, including the case of Schwarzschild potential, we have
that
$$
\int_0^T \| \langle A \rangle ^{1/2}  F_M^{\pm} \langle x \rangle ^{-2} u \| ^2 dt \leq CE(u).
$$
\endproclaim
\demo{Proof}
$$
$$
The proof for the $u_l$ part is completed by the above theorem, on
noticing that for $\sigma =2$, we have that
$$
\int_0^T\| \tilde g(A) \langle x \rangle ^{-2} u_{_l} \| ds \leq
C\int_0^T T^{-1/2} E^{1/2}(u_{_l}) dt \leq CT^{1/2} E(u_{_l}).
$$
\enddemo

\subhead Analytic Repulsiveness of the Schwarzschild potentials
\endsubhead

When the potential vanishes at - infinity, exponentially fast, the
situation is complicated by the fact that, even though
$$
-x \cdot \nabla V \geq 0 \text{ at infinity,}
$$
in general, $V_\beta$ is not positive, but oscillates no matter how
small $\beta$ is:

For $V(x) = e^{-x}$ for $x >> M,$
$$
V_\beta = 2 Im e^{-e^{-i\beta}x} = 2 Im e^{- x \cos\beta} e^{+ xi \sin
\beta}
$$
$= 2 e^{-x\cos \beta} (+\sin(x \sin\beta))$ which decays exponentially, but
oscillates with period $(\sin \beta)^{-1}$.

So, to prove analytic repulsiveness, we need to show that
$$
2\sin 2 \beta p^2 + V_\beta \sim 2\sin2\beta p^2 + 2 e^{-x\cos\beta} \sin( x
\sin \beta)
$$
is a positive operator.

\proclaim{ Theorem 7.10}

Suppose $V(x)$ is repulsive: $-x \frac{\partial V}{\partial x} \geq
f^2 (x) > 0,$ one hump potential, with non degenerate maximum.

Suppose, moreover, that $V_\beta$ exists and is analytic for all
$|\beta|$ sufficiently small, and
$$
|V_\beta(x)| \leq C e^{-\delta x}, x > x_0, \text{ for some } x_0 >
 0.\tag i
$$
$$
|V_\beta(x) | \leq C\langle x\rangle^{-2-a} \text{ for all } |x| > +
  x_0, \text{ some } a > 0.\tag ii
$$
condition (ii) can be replaced by condition iii):
$$
V_\beta(x) \geq f^2(x, \beta) > 0 \text{ for } x < - x_0.\tag iii
$$
\endproclaim

Then, $V$ is analytic-repulsive, and
$$
2\sin 2 \beta p^2 + V_\beta \geq
\delta_0 \langle x\rangle^{-2}.\tag 7.7
$$

\noindent{\bf Remarks}

The condition on $V$ implicitly implies that $V$ has a (dilation)
analytic extension from $\Bbb R$ to the domain
$$\{e^{i\beta'} x\big| |\beta'| \leq \beta, x-\text{real}\}.
$$

\demo{Proof}

Using the fundamental theorem of calculus and Taylor series expansion,
we write $V_\beta$ as
$$
\aligned
& V_\beta(x) \int^\beta_{-\beta} 2\{ Im \frac{\partial}{\partial s}
V(e^{-is} x)\} ds\\
& =-2 x V'(x)\beta + 2\int^\beta_{-\beta} Re \{e^{-is'-is} V''(e^{-is'} x)\} s
x^2 ds\\
& \geq -2 xV'(x) \beta - \beta^2\frac{|x|^2}{2}\sup_{|s'|\leq\beta}
|V''(e^{-s'} x)|.\endaligned
\tag 7.8
$$
Using the following Cauchy estimates
$$
\vert f^{(n)}(z_0)\vert\leq \frac{n!}{R^n}\sup_{\vert z-z_0\vert=R} \vert f(z)\vert,
$$

 we have that:
$$
V_\beta\geq -2 xV'(x)\beta-\beta^2 |x|^2 2|\sup_{s'}\sup_{|z - e^{is'}
x |=1} |V(z)|
$$
where $-\beta \leq s' \leq \beta$.

Using condition (i) of the theorem, it follows that for large $x$
positive, $x> x_0$:
$$
V_\beta(x) \geq -2 x V'(x)\beta - c\beta^2|x|^2 e^{-\delta x}.\tag 7.9
$$
For $x< x_0$:

Since $V(x)$ is assumed to be a one hump potential, $x V'(x)$ is
strictly positive away from zero, and (non degenerate case)
$$
xV'(x) \sim \frac{1}{2} a^2 x^2 \text{ near zero}.
$$
Here $x=0$ is the top of the hump of $V(x)$.  Since $|V(z)|\leq
C\langle x \rangle^{-2}$, choosing $\beta$ sufficiently small, we have
that for all $|x|\leq x_0$.
$$
-2x V'(x)\beta - 2\beta^2 |x|^2 \sup_{|s'|\leq \beta}\sup_{|z-e^{is'}
 x|=1} |V(z) | \geq - xV'(x) \beta.\tag 7.10
$$
For $x$ large, negative, we use the Cauchy estimate with $
|z-e^{is'} x | =1$ replaced by a  circle, which encloses $e^{is'} x$
of radius $\sim \delta(\beta|x|)^{1-\eta};  \eta, \delta$ small, so that
$z$ is in the domain of analyticity.

Then, we have that for $x< - x_0$:
$$
\aligned
V_\beta (x) & \geq - 2 x V'(x) \beta - c\beta^{2\eta} |x|^{2\eta} \langle x
\rangle^{-2} \langle x \rangle^{-2-a} /\delta^2\\
            & \geq -2 x V'(x) \beta - c\beta^{2a} \langle x\rangle^{-2-\varepsilon}/\delta^2
\endaligned \tag 7.11
$$
for$ 2\eta < a$:

 If condition (iii) is satisfied then
$$
V_\beta(x) > 0 \text{ for all } x < - x_0.
$$

Now, since $xV'(x)> 0$ for $x \neq 0$, by choosing $\beta$
sufficiently small, we have that
$$
V_\beta(x) \geq \frac{\beta}{2} f^2 (x) \text{ for all } x <x_0 \tag
7.12
$$
and
$$
|V_\beta(x)| \leq C e^{-\delta x} ,  x <x_0.\tag 7.13
$$

So, to complete the proof, we need to show that $2 \sin 2 \beta p^2 +
V_\beta (x) > 0$.

To this and, we use the uncertainty principle, which, in one
dimension, gives (in one of its forms... $ [BSt]$)
$$
p^2 + \lambda \chi_I \geq \frac{C(\lambda, I)}{1+|x|^2} \tag 7.14
$$
where $\chi_I$ is the characteristic function of the interval
$ I:$

By 7.9 - 7.14,
$$
2 \sin 2 \beta p^2 + V_\beta(x) \geq 2\sin 2 \beta \left(p^2 + \frac{1}{2}
f^2 (x) \chi (x < x_0) \right) - C\chi(x> x_0)e^{-\delta x} \beta^2 x^2
$$
$$
\geq 2 \sin\beta\left(p^2 + \frac{1}{2} f^2(x) \chi(x< x_0)\right) - C\beta ^2
\chi(x> x_0) e^{-\delta x} x^2(1+x^2) \frac{1}{1+x^2}
$$
$$
\geq\frac{1}{1+x^2}
\left[
\beta C (x_0, f^2) - \chi (x> x_0) C\beta^2
e^{-\delta x_0/2} e^{-\delta(x-x_0/2)} x^2 (1+x^2)
\right]
$$
$$
\geq \frac{\beta}{2} C(x_0, f^2)(1+x^2)^{-1},
$$
by choosing $\beta$ small, and by choosing $x_0$ large enough so that
$x_0>> \frac{2}{\delta}$, to get
$$
C\beta^2 e^{-\frac{\delta x_0}{2}} \delta^{-4} < \frac{\beta}{2} C(x_0,
f^2)
$$
which is possible, since increasing $x_0$ only increases the value of
$C(x_0, f)$.
\enddemo
 \qed

\proclaim{ Theorem 7.11} (Improved Local Decay for Schwarzschild
potentials)

Let
$$
H=-\Delta + V_\ell (x) \quad  x \in \Bbb R.
$$
$$
V_\ell \text{ is defined in } (7.1) - (7.3) (x\equiv r_*).
$$
Then, the following local decay estimate holds:
$$
\ell \geq 1: \int^T_0 dt \|J u (x, t)\|^2 \leq C E^{1/2} (u)
E^{1/2}\left( \langle p \rangle^{-1}pu\right)\tag i
$$
$$
J=J(x) , |J(x)| \leq (1+ x^2)^{-1/2-\delta}.
$$
$$
\aligned
\ell = 0; & \\
&\int^T_0 \|J u (x, t) \|^2 dt \leq C E^{1/2} (u) E^{1/2} \left(\langle p
\rangle^{-1} pu \right) \\
\text{ with } & \\
 &J=J(x), |J(x)|\leq (1+x^2)^{-3/2 - \delta}.
\endaligned
\tag ii
$$
\endproclaim

\demo{Proof}

For each $\ell$, the potential $V_\ell(x)$ is a one hump function
[B-Sof1], and has analytic continuation [Bac-Bac, Zw] for all $\beta$
sufficiently small.

Moreover,it satisfies the conditions of Theorem 7.2 [Bac-Bac, Zw], and
$-x\nabla V\geq f^2(x) $ is also known [B-Sof1].

So, applying Theorems 7.2, 6.4 and proposition 7.1 the result follows.
\qed
\enddemo

{\bf Example } {\it Negative Potentials in 3 dimensions}
$$
V(x) = - \left(\frac{a}{b^2 + x^2}\right)^2 \text{ in three
dimensions}.
$$

Then
$$
-(2-\varepsilon) V- x\cdot \Delta V=
\frac{-(2-\varepsilon) (b^2 + r^2) + 4r^2}  {b^2 + r^2}
V = \frac{ (2+\varepsilon) r^2-b^2(2-\varepsilon)}{b^2 + r^2} V
$$
since
$$
-x\cdot \nabla V = \frac{+ 4 r^2}{(b^2 + r^2)} V.
$$

The above expression is negative for $r^2\geq b^2
\frac{2-\varepsilon}{2+\varepsilon}$.

Hence
$$
\aligned
& 2 p^2 - x\cdot \nabla V = (2-\varepsilon)(p^2+V) + [-x \cdot \nabla
V - (2 - \varepsilon) V ] + \varepsilon p^2\\
&
\geq (2-\varepsilon) \delta |x|^{-2} + \frac{\varepsilon}{4}|x|^{-2} -
\frac{(2+\varepsilon) r^2 - b(2-\varepsilon)}{b^2 +r^2}
\frac{a^2}{(b^2 + r^2)^2} \\
& >0, \text{ for } (a/b) \text{ sufficiently small}.
\endaligned
$$

{\bf Example }  -{\it Addition of Humps}

This example is typical to the problem of constructing a propagation
observable with no $\ell$ dependence for the Schwarzschild/Kerr
problem, for example.

Here, I consider a simple example, leaving the general case to other
works.

So, let
$$
V(x) = \frac{2}{1+|x|^2} + a\frac{1}{1+ |x-b|^3}
$$
for $a, b > 0, -\infty < x < \infty$.

Then
$$
-x\cdot \nabla V = \frac{4 r^2}{(1+ r^2)^2} +
 \frac{3a|x-b|^3}{(1+|x-b|^3)^2}
+ \frac{3 a b \,  sgn(x-b) |x- b|^2}{(1+| x-b|^3)^2}.
$$

This expression  may be negative for $0\leq x \leq b$.  It is
negative near $x = 0 , x > 0$ since the last term dominates.

However, using the localized uncertainty principle, Lemma 4.3b, it
follows that
$$
\frac{1}{4} p^2 + \frac{4r^2}{(1 + r^2)^2} \geq \frac{1}{4}
\frac{1}{(1+ r^2)^2}.
$$
Therefore, we can easily arrange
$$
2 p^2 -x\cdot \nabla V \geq \varepsilon \langle x \rangle^{-4}
$$
by choosing $a$ small or $b$ small.

\noindent{\bf Other Perturbations}
\medskip

All the previous examples will still satisfy the local decay estimates
under the addition of a small, fast decaying, possibly time dependent
perturbation, $W(x, t)$, provided $W_\beta(x, t)$ is well defined for
small $\beta$, and satisfies the same size and decay conditions.

\medskip

\head Section 8\endhead
\subhead
8. High Angular Momentum Bounds\endsubhead
\medskip
In this section we demonstrate an application to Schwarzschild scattering, for large angular momentum.
It is by no means supposed to be comprehensive, and the genral results, including pointwise estimates will be developed elsewhere.
In the previous sections we did not follow the dependence of the decay estimates on the $\ell$ dependence.
Here, we will consider the angular dependence of the previously obtained decay estimates, for the Schwarzschild potential and for the case where
$$
V_{\ell}=\ell^2V(x),
$$
with $V(x)$ analytic repulsive. This is motivated by the case of extreme Reissner Nordstrom Blackhole manifold.
Our main goal is to show, that for large $\ell,$ the local decay estimate holds, with a factor of $\ell,$ up to log correction. Previously, this was proved in [B-Sof3,4 ], by a complicated generalized phase-space analysis.
We begin with the following preliminary results, that follow directly from applying the previous estimates.
First, we note that, in the Schwarzschild case, the behavior of the potential at large negative $x_*,$ is $\ell^2$ times an exponentially decaying function.
Therefore, to insure that such a potential is repulsive analytic, we need to choose $\beta,$ in the definition of the PROB, to be smaller than $c(\ln \ell)^{-1},$ for some sufficiently large positive $c.$
Then, we have the following estimates:

\proclaim{Proposition 8.1}
Let $$H=-\Delta +\ell^2 V(x),\tag8.1$$
with $V(x)$ analytic repulsive, for $\beta \leq \beta_0(\ell).$
Then, we have the following PRES:
$$
\int_0^T \| H_{\beta}^{1/2}\tilde g(A/R)u\|^2 dt \leq cE(u)^{1/2}E^{1/2}(\left <p\right>^{-2}pu),\tag8.2
$$
$$
\int_0^T \| Q_{\beta}^{1/2}\tilde g(A/R)\dot u\|^2 dt \leq cE(u),\tag8.3
$$
where we define $Q:=\sqrt H.$
$$
\int_0^T \| J(x)\ell u\|^2 dt \leq cE(u),\tag8.4
$$
where $ J(x)=cx \left <x\right>^{-2}(1+e^{bx})^{-1},$
with $b$ positive.
\endproclaim
\medskip

The proof of the above statements follows from application of the previous PRES to the hamiltonian defined in (8.1).
\smallskip

{\bf Sketch of Proof}
The estimate (8.2) follows by using the PROB $\tanh (A/R)$, together with Theorem (4.4). $R$ is chosen large enough, depending on $\ell$, to insure the positivity of $ H_{\beta}=\frac{i}{2}\left[H^{[-\beta]} - H^{[\beta]}\right].$ Here we note that the algebraic proof of Theorem (4.4) applies verbatim with $H$ replacing $V$.
The resulting PRES is the estimate (8.2).

The estimate (8.3) follows by repeating the above argument for the Schr\"odinger type equation, with hamiltonian given by $Q=\sqrt H.$ To this end, we note that the equation satisfied by the function $\dot u$, is given by
$$
\dot u(x,t)=\sin(Qt) (Qf)+\cos(Qt) g,
$$
with $Qf, g$ in $L^2$.
 The sine and cosine functions are linear combinations of $e^{\pm iQt},$ which is the propagator of the Schr\"odinger equation with hamiltonian
$\mp Q.$
Applying as above Theorem (4.4) and the resulting PRES, we obtain (8.3).
To prove the estimate (8.4), we write the PRES for the following PROB
$$
G:= b(x)\partial_x +\partial_x b(x),
$$
with $b(x)=x/\left <x\right >.$
Then,  we obtain a positive term from the commutator with the potential part, of the form
$$
-2b(x)\ell^2 V'(x),
$$
together with two terms from the commutator with the Laplacian part of the hamiltonian.
One term is positive, and is of second order in the radial derivative; the other is localized in x. This localized term, has coefficient of order 1, that is, independent of $\ell.$ It comes from $b'''$ term in the commutator.
Since we proved that for such localized weight function the PRES holds, the result follows.
\medskip

\proclaim{Proposition 8.2}
Under the same assumptions of Theorem (8.1), we have the following PRES:($r_0>0$)
$$
\int_0^T \left \{\|xF(|x|\leq r_0)\ell \tilde g(A/R)u\|^2+\|F(|x|\leq r_0)\sqrt\ell \tilde g(A/R)u\|^2\right \}\leq
c[\left<F(A/R)\right >_{u}^{\text{Heis}}-\left<F(A/R)\right >_{u_0}^{\text{Heis}}].\tag8.5
$$

\endproclaim
The above proposition is a consequence of previous decay estimates, with $V$ replaced by $\ell V.$
The second term on the rhs, is bounded, with a loss $\ell^{1/2},$ to eliminate the vanishing $x$ factor.
This follows from application of the uncertainty principle, as in [DSS2].
\medskip

We will be able to get the desired estimate from this last bound, on using it with $u\rightarrow Q^{1/2}u,$
and using the fact that $Q>F(|x|\leq r_0)\ell.$
The resulting estimate is restricted to the support of the operator $g(A/R).$
To this we show how to remove this projection from the estimate.
\medskip

\proclaim {Proposition 8.3}
Under the same assumptions of Theorem (8.1), we have the following PRES:($r_0>0$)
$$
\int_0^T \|\left <x\right >^{-a} |A|^{1/2}F_M(A/R)u\|^2 dt\leq c\Re[\left< \left <x\right >^{-a'}F(A/R)\right >_{u}^{\text{Heis}}-\left< \left <x\right >^{-a'}F_M(A/R)\right >_{u_0}^{\text{Heis}}],\tag8.6
$$
where $a=(a'+1)/2.$

\endproclaim

This follows, as before, by using the following PROB, similar to the one used before, with similar computations:
$$
\left<x\right>^{-a'}F_M +F_M\left<x\right>^{-a'}.
$$
\medskip
We can then use the PRES to remove the cutoff function $g(A/R)$ from the PRES of proposition (8.2), except that we need to bound this error term by a quantity that is of order $\ell^{-2},$ up to possibly log corrections, for large $\ell.$
The first power of $\ell$ comes from , as before, by applying the above proposition to $u\rightarrow Q^{1/2}u.$
To obtain another power of $\ell$, we use the {\bf redeeming property} of the Heisnberg identity for the wave equation: If $N$ is a symmetric reality preserving linear operator, then:
$$
\left< N\right >_{u}^{\text{Heis}}=0.\tag 8.7
$$
In particular, this holds for $N=1, f(x),g(|p|),fg+gf$ with $f,g$ real valued functions.
Therefore,
 as noted before, we have
 $$
 \aligned
 &\left<F_M(A/R)\right >_{u}^{\text{Heis}}=\\&
 \left<F((A-M)/R\sim K_0)\right >_{u}^{\text{Heis}}\\&
 +\left<F(|(A-M)/R|\ge K_0)F_M\right >_{u}^{\text{Heis}}.
 \endaligned
 $$
 $$
 \left<F(|(A-M)/R|\ge K_0)\right >_{u}^{\text{Heis}}=\left<(F_1(|x|\leq C)+\bar F_1(|x|\ge C))F(|(A-M)/R|\ge K_0)\right >_{u}^{\text{Heis}},
 $$
 $$
 \aligned
 &
 \left<\bar F_1(|x|\ge C)F(|(A-M)/R|\ge K_0)\right >_{u}^{\text{Heis}}\\&=
 \left< F_x(F_p+\bar F_p)F(A)\right>_{u}^{\text{Heis}}=\left< F_xF_p\right>_{u}^{\text{Heis}}+\left< F_x \bar F_p\right >_{u}^{\text{Heis}}+O(\ell^{-2})=O(\ell^{-2})+\left< F_x \bar F_p\right >_{u}^{\text{Heis}}.
 \endaligned
 $$
 Here, $F_x\equiv \bar F_1(|x|\ge C),$ $F_p\equiv F_p(|p|\ge \delta \ell),$ $\delta\ge0.$
 $F(A)\equiv F(|(A-M)/R|\ge K_0).$
 We are therefore left with controlling (by $O(\ell^{-2})$) the regions of phase space:
 $$
 \aligned
 & \left < F_1(|x|\leq C)F(|A-M|\ge K_0)\right >_{u}^{\text{Heis}}\\
 & \left < F(|A-M|\leq K_0)F_M\right >_{u}^{\text{Heis}}\\
 & \left< F_x \bar F_p F(|A-M|\ge K_0) \right >_{u}^{\text{Heis}}.
 \endaligned \tag 8.8
 $$
 To complete the proof of the main estimate with $O(\ell^{-2})$ decay, up to logarithmic corrections in $\ell,$
 we need to bound the above three terms of the formula (8.8), which are referred below as terms I,II,III, with $u\rightarrow Q^{1/2}u,$
by $\ell^{-1}E(u).$

To this end, we estimate the scalar product as follows:
$$
|\left< Q^{1/2}u,G_1 G_2Q^{1/2}\dot u\right >| \leq C \|G_2\dot u\| \|G'_1 Qu\| ,\tag 8.9
$$
for  generic operators $G$, and with $G'_1\equiv G_1+[Q^{1/2},G_1]Q^{-1/2}.$

Estimate of I:
$$
\left <F_1(|x|\leq C)F(|A-M|\ge K_0)\right >=\left<F_1\right>-\left<F_1\bar F\right>=-\left<F_1\bar F\right>.\tag 8.10
$$
Therefore,  this last term is bounded by $O(\ell^{-1}\ln \ell)E(u),$ (with $u\rightarrow Q^{1/2}u,$), by applying Proposition 8.2.

Estimate of II:
$$
\left<FF_M\right>=\left<FF_MF_p\right>+\left<FF_M\bar F_p F(|x|\leq C)\right>+O(\ell^{-2}),\tag 8.11
$$
since, as we will show below, we only need to consider initial data with $F(1/2 \ell \leq Q\leq 2\ell)u=u,$
and the localization lemmas below, that imply $\bar F_p F(|x|\ge C)F(Q\ge (1/2) \ell)=O(\ell^{-2}).$
 The second term, on the right hand side of equation (8.11), is bounded by $O(\ell^{-1}\ln \ell)E(u),$ (with $u\rightarrow Q^{1/2}u,$), by applying Proposition 8.2 as before.

 The first term,on the right hand side of equation (8.11),is bounded by $O(\ell^{-1}\ln \ell)E(u),$ (with $u\rightarrow Q^{1/2}u,$), by applying Proposition 8.1,since $F_pQF_p\ge c\ell F_p.$

 Estimate of III:
 $$
 \left< F_x(|x|\ge C)F(|p|\leq \delta \ell)\right>=O(\ell^{-2}),
 $$
 again, by the localization lemmas below.

 \proclaim{Lemma 8.4} {\bf Localization lemma}
 Let $H=-\partial_x^2+\ell^2V(x)$ be defined as before in this section.
 Furthermore, we normalize $V(0)=1.$
 Then,for all $n>0,$
 $$
 F(H\ge (1/2)\ell^2)F(|x|\ge c)F(|p|\leq \delta \ell)=O(\ell^{-n}),\tag i
 $$
 for all $c$ large enough. $ \delta < 1/2.$
 $$
 F(H\ge2\ell^2)F(|p|\leq \delta\ell)=O(\ell^{-n}),\tag ii
 $$
 $$
 F(H\leq (1/2)\ell^2)F(|x|\leq \delta)=O(\ell^{-n}),\tag iii
 $$
 \endproclaim
 \demo{Proof}
 The proof follows the method of proving the Localization Lemma of [Sig-Sof1,2]:
 i) Let us denote by $g=F(|x|\ge c)F(|p|\leq \delta \ell).$
 Then, we define
 $$
 \bar H \equiv g^* Hg \tag 8.12
 $$
 We have:
 $$
 \bar H =F_pF_x(p^2+V(x)\ell^2)F_pF_x \leq \delta ' \ell^2.
 $$
 Therefore $ F(\bar H\ge (1/2)\ell^2)=0.$
 Then, with $\tilde g g=g,$ and all positive integers $k,k',$
 $$
 \aligned
 &\tilde g^* F(H\ge (1/2)\ell^2)\tilde g=\tilde g^*\{ F(H\ge (1/2)\ell^2)- F(\bar H\ge (1/2)\ell^2)\}\tilde g\\
 &=\tilde g^*\int\hat F(\lambda)e^{iH\lambda}\int_0^{\lambda}e^{-iHs}[H-\bar H]e^{is\bar H}d\lambda ds \tilde g\\
 &=\tilde g^*\int\hat F(\lambda)\int_0^{\lambda}Ad_g^{(k)}(e^{iH(\lambda-s)}[H-\bar H]Ad_g^{(k')}[e^{is\bar H}d\lambda] ds \tilde g.
 \endaligned \tag 8.13
 $$
 Direct computation shows that $[g,H],[g,\bar H] =O(\ell).$
 Therefore,
 $$
 [g,e^{itH}]=ce^{itH}\int_0^te^{-isH}O(\ell)e^{+isH}ds=O(t\ell).\tag 8.14
 $$
 By repeatedly commuting $g$ through the above expression, and using the fact that $[g,O(\ell)]=O(1),$
 it follows that the multicommutators in equation 8.13
 are bounded by $c_kt^k \ell^k,$, for some constants $c_k$, depending only on the sharpness of the functions defining $g.$
 Since $H-\bar H=O(\ell^2),$ direct estimate of the $L^2$ norm of the rhs of equation 8.13 gives:
 $$
 \aligned
 &\|\tilde g^* F(H\ge (1/2)\ell^2)\tilde g\|=\|\tilde g^*\{ F(H\ge (1/2)\ell^2)- F(\bar H\ge (1/2)\ell^2)\}\tilde g\|\\
 &\leq c_n \int |\lambda^n \hat F(\lambda)|d\lambda O(\ell^{2+n-1}).
 \endaligned \tag 8.15
 $$
 Finally, using the construction of the function $F,$
 we have that
 $$
 \int |\lambda^n \hat F(\lambda)|d\lambda \leq c_n\ell^{-2n+2}. \tag 8.16
 $$
 \enddemo
Putting it all together, we establish the following improved local decay estimate, for large $\ell:$

\proclaim{Theorem 8.5}
For the Hamiltonian with the Schwarzschild potential $\ell^2 V(x),$ with $V$ analytic repulsive, we have the following estimate:
$$
\int_0^T \|F(|x|\leq r_0)\ell u\|^2 \leq c\ln \ell E(u). \tag 8.17
$$
\endproclaim

\noindent{\bf Acknowledgments}

This work is partially supported by NSF grant DMS-0903651

\enddocument
\ref
\key
\by
\paper
\jour
\vol
\yr
\pages
\endref
\medskip
\ref
\key
\by
\paper
\jour
\vol
\yr
\pages
\endref
\medskip

\ref
\key
\by
\paper
\jour
\vol
\yr
\pages
\endref
\medskip
\end

\enddocument
\end

enddemo\enddocument

\enddocument

\end